\documentclass{svmult}
\usepackage{makeidx}         % allows index generation
\usepackage{graphicx}        % standard LaTeX graphics tool
\usepackage{multicol}        % used for the two-column index
\usepackage[bottom]{footmisc}% places footnotes at page bottom
\makeindex             % used for the subject index
\newtheorem{thm}{Theorem}%[section]
\newtheorem{lem}[thm]{Lemma}
\newtheorem{prop}[thm]{Proposition}
\newtheorem{cor}[thm]{Corollary}
\newtheorem{rmk}[thm]{Remark}

\begin{document}

\baselineskip14pt

\title*{Variation and Rough Path Properties of Local Times of L\'evy Processes}
\titlerunning {Local time rough path}
\author{Chunrong Feng\inst{1},
Huaizhong Zhao\inst{2}}
\authorrunning{C.R. Feng and H.Z. Zhao}
\institute{Department of Mathematics, Shanghai Jiao Tong University,
Shanghai, 200240, China \and Department of Mathematical Sciences,
Loughborough University, LE11 3TU, UK \\ \texttt{fcr@sjtu.edu.cn},
\texttt{H.Zhao@lboro.ac.uk}} \maketitle
%\footnotetext{This is the
%draft version with more details of the paper submitted to COMPTES
%RENDUS MATHEMATIQUE}
\begin{abstract}

In this paper, we will prove that the local time of a L\'evy process
is of finite $p$-variation in the space variable in the classical
sense, a.s. for any $p>2$, $t\geq 0$, if the L\'evy measure
satisfies $\int_{R\setminus \{0\}}(|y|^{3\over 2}\wedge
1)n(dy)<\infty$, and is a rough path of roughness $p$ a.s. for any
$2<p<3$ under a slightly stronger condition for the L\'evy measure.
Then for any function $g$ of finite $q$-variation ($1\leq q <3$), we
establish the integral $\int _{-\infty}^{\infty}g(x)dL_t^x$ as a
Young integral when $1\leq q<2$ and a Lyons' rough path integral
when $2\leq q<3$. We therefore apply these path integrals to extend
the Tanaka-Meyer formula for a continuous function $f$ if $\nabla
^-f$ exists and is of finite $q$-variation when $1\leq q<3$, for
both continuous semi-martingales and a class of L\'evy processes.

\vskip10pt

\end{abstract}
\vskip-10pt

\newcounter{bean}
 \renewcommand{\theequation}{\arabic{equation}}

\section {Introduction}

The variation of a stochastic process is a classical problem of
fundamental importance in probability theory. There are two kinds of
variations, namely in the sense of probability and in the classical
sense.  The quadratic variation of a Brownian motion (a martingale)
in the sense of probability made it possible to define It${\hat {\rm
o}}$'s stochastic integral of a square integrable progressive
process with respect to the Brownian motion (the martingale) (L\'evy
\cite{levy}, It${\hat {\rm o}}$ \cite{ito1}, Kunita and Watanabe
\cite{kunwat}). The classical $p$-variation of the Brownian motion
and its L\'evy area, $p>2$, led to Lyons' pathwise approach to
stochastic differential equations (Lyons \cite{terry2},
\cite{terry1}). Local time is an important and useful stochastic
process. The investigation of its variation and integration has
attracted attentions of many mathematicians. Similar to the case of
the Brownian motion, the variation of the local time of a
semimartingale in the spatial variable is also fundamental in the
construction of an integral with respect to the local time. There
have been many works on the quadratic or $p$-variations (in the case
of stable processes) of local times in the sense of probability.
Bouleau and Yor (\cite{bou}), Perkins (\cite{perkins}), first proved
that, for the Brownian local time, and a sequence of partitions
$\{D_n\}$ of an interval $[a,b]$, with the mesh $|D_n|\to 0$ when
$n\to \infty$, $\lim\limits_{n\to \infty} \sum \limits_{D_n}
(L_t^{x_{i+1}}-L_t^{x_i})^2=4\int _a^bL_t^xdx$ in probability.
Subsequently, the process $x\to L_t^x$ can be regarded as a
semimartingale (with appropriate filtration). This result allowed one
to construct various stochastic integrals of the Brownian local time
in the spatial variable. See also Rogers and Walsh \cite{rog}.
Numerous important extensions on the variations, stochastic
integrations of local times and It$\hat {\rm o}$'s formula have been
done, e.g. Marcus and Rosen \cite{rosen1}, \cite{rosen}, Eisenbaum
\cite{eisenbaum1},
 \cite{eisenbaum11},
Eisenbaum and Kyprianou \cite{eisenbaum12},
 Flandoli, Rosso and Wolf \cite{frw}, F\"ollmer,
Protter and Shiryayev \cite{Protter}, F\"ollmer and Protter
\cite{Protter2}, Moret and Nualart \cite{nualart}. In their
extensions of It${\hat {\rm o}}$'s formula, the integrals of the
local time are given as stochastic integrals in nature, for example
as forward and backward stochastic integrals. What had been missing
here in the literature mentioned above was the variation of the
local time in the classical sense and whether or not the local time
is a rough path (its meaning will be made precise later). In Feng
and Zhao \cite{Zhao33}, we proved that for a
 continuous semi-martingale,
 its local time $L_t^x$ is of finite
$p$-variation in the classical sense in $x$  for any $t\geq 0$,  a.s. for any $p>2$, i.e.
 \begin{eqnarray}
 \sup_{D}\sum _{i=0}^{n-1}|L_t^{x_{i+1}}-L_t^{x_i}|^p<\infty,
 \end{eqnarray}
 where the supremum is taken over all finite partition $D_{(-\infty,+\infty)}=\{-\infty<x_0<x_1<x_2<\cdots <x_n<\infty\}$.
  This allowed us to define the path integral $\int _{-\infty}^\infty g(x)d_xL_t^x$
 as a Young integral, for any $g$ being of a finite $q$-variation for a number $q\in [1,2)$. The main purpose of this paper is to solve the rough path part of the problem,  for the local times
 of both continuous semi-martingales and a class of L\'evy processes.

 As a first step, we consider the classical $p$-variation for the local time of the L\'evy process which is represented by the following
L\'evy-It${\hat {\rm o}}$ decomposition
\begin{eqnarray}\label{tanakaf1}
X_t&=&X_0+\sigma B_t+bt+\int_0^{t+}\int_{R\setminus \{0\}}y 1_{\{|y|\geq
1\}} N_p(dsdy)+\int_0^{t+}\int_{R\setminus \{0\}}y 1_{\{|y|< 1\}}
\tilde N_p(dsdy).
\end{eqnarray}
 This is a non-trivial problem as the jumps, especially the small jumps, create
a lot of difficulties in estimating the increment
$L_t^{x_{i+1}}-L_t^{x_i}$. Recall that for a general semimartingale
$X_t$,  $L=\{L_t^x; x\in R\}$  is defined from the following formula
(Meyer \cite{meyer}):
\begin{eqnarray}\label{tanakaf2}
\int_0^t g(X_s)d[X,X]^c_s=\int_{-\infty}^\infty g(x)L_t^x dx,
\end{eqnarray}
where $[X,X]_.^c$ is the continuous part of the quadratic variation process $[X,X]_.$. There is a different
notion of local time defined as the Radon-Nikodym derivative of the
occupation measure of $X$ with respect to the Lebesgue measure on
$R$ i.e.
\begin{eqnarray}\label{tanaka3}
\int_0^t g(X_s)ds=\int_{-\infty}^\infty g(x)\gamma_t^x dx,
\end{eqnarray}
for every Borel function $g: R\to R^+$. For the L\'evy process
(\ref{tanakaf1}), if $\sigma\ne 0$, $L_t^x$ and $\gamma _t^x$ are the same
(up to a multiple of a constant). In case $\sigma=0$ e.g. for a stable
process,  there is no diffusion part so these two definitions
are different. In fact, in this case $L_t^x=0$. The increment of
 $\gamma_t^.$ for stable processes was considered by
Boylan \cite{bo}, Getoor and Kesten \cite{geke} and Barlow
\cite{barlow}, using potential theory approach, in order to
establish the continuity of the local time in the space variable.
The first main task of this paper is to prove that when $\sigma \ne
0$, and (\ref{levy110}) is satisfied, for any $p>2$ and $t\geq 0$,
the process $x\mapsto L_t^x$, is of finite $p$-variation in the
classical sense almost surely. Both our result and our method are new in literature.
Here the Tanaka formula is directly used in our approach.
As a direct application, one can define the path integral $\int
_{-\infty}^{+\infty} g(x)dL_t^x$ as a Young integral for any $g$
being of bounded $q$-variation for a $q\in [1,2)$. It is noted when
$q\geq 2$, Young's condition
 ${1\over p}+{1\over q}>1$ is broken.

The main task of  this paper is to construct a geometric rough path
over the processes $Z(x)=(L_t^x,g(x))$, for a deterministic function
$g$ being of finite $q$-variation when $q\in [2,3)$. This implies
establishing the path integrals $\int_{-\infty}^\infty L_t^xd_xL^x_t
$ and $\int_{-\infty}^\infty g(x)d_xL^x_t $. For these two
integrals, all the classical integration theories of Riemann,
Lebesgue and Young fail to work. To overcome the difficulty, we use
the rough path theory pioneered by Lyons,
 see \cite{terry2}, \cite{terry1}, \cite{terry}, also \cite{lejay}.
 However, our $p$-variation result of the local time does not automatically make the desired rough path exist or
 the integral well defined, though it is a crucial step to study first. Actually further hard analysis is needed
 to establish an iterated path integration theory for $Z_.$.
 First we introduce a piecewise curve of bounded variation as a generalized Wong-Zakai approximation to the stochastic process $Z_.$. Then
 we define a smooth rough path by defining the iterated integrals of the piecewise bounded variation process. We need to prove the smooth
 rough path converges to a geometric rough path ${\bf {Z}}=(1,{\bf { Z}}^1,{\bf { Z}}^2)$ when $1\leq q<3$.  For this,
an important step is to compute
$E(L_t^{x_{i+1}}-L_t^{x_i})(L_t^{x_{j+1}}-L_t^{x_j})$, and obtain
the correct order in terms of the increments $x_{i+1}-x_i$ and
$x_{j+1}-x_j$, especially in
 disjoint intervals $[x_i, x_{i+1}]$ and $[x_j, x_{j+1}]$ when $i\ne j$.
 Actually, this is a very challenging task. In this analysis, one can see that
 the main difficulty is from dealing with jumps, especially the small jumps of the
 process. One can also see that (\ref {levy110}) is not adequate to
 construct the geometric rough path, a slightly stronger condition
 (\ref {levynew36}) is needed here.
 Using this key estimate, we can establish the geometric rough path
   ${\bf {Z}}=(1,{\bf { Z}}^1,{\bf { Z}}^2)$.
Then from Chen's identity, we define the following two integrals
 \begin{eqnarray}\label{tanakaf3}
\int_a^b L_t^xd L_t^x=\lim_{m(D_{[a,b]})\to 0}\sum_{i=0}^{r-1}(({\bf
Z}_{x_i,x_{i+1}}^2)_{1,1}+ L_t(x_i)(L_t^{x_{i+1}}- L_t^{x_i}))
\end{eqnarray}
and
\begin{eqnarray}\label{tanakaf4}
\int_a^b g(x)d L_t^x=\lim_{m(D_{[a,b]})\to 0}\sum_{i=0}^{r-1}(({\bf
Z}_{x_i,x_{i+1}}^2)_{2,1}+g(x_i)( L_t^{x_{i+1}}- L_t^{x_i})).
\end{eqnarray}
 Note that the
Riemann sum $\sum\limits _{i=0}^{r-1}L_t^{x_i}( L_t^{x_{i+1}}-
L_t^{x_i})$ and $\sum\limits _{i=0}^{r-1}g(x_i)( L_t^{x_{i+1}}-
L_t^{x_i})$ themselves may not have limits as the mesh $m(D_{[a,b]})
\to 0$. At least there are no integration theories, rather than
Lyons' rough path theory, to guarantee the convergence of the
Riemann sums for almost all $\omega$. Here it is essential to add
L\'evy areas to the Riemann sum. Furthermore, we can prove if a
sequence of smooth functions $g_j\to g$ as $j\to \infty$, then the
Riemann integral $\int_a^b g_j(x)dL_t^x$ converges to the rough path
integral $\int_a^b g(x)dL_t^x$ defined in (\ref{tanakaf4}). It is
also noted that to establish (\ref{tanakaf3}), one only needs
(\ref{levynew37}). This is true as long as the power of $|y|$ in the
condition of L\'evy measure is anything less (better) than ${3\over
2}$. It seems to us ${3\over 2}$ is a critical value here. Our
conjecture is that when the L\'evy measure satisfies
$\int_{R\setminus \{0\}} (|y|^{3\over 2} \wedge 1)n(dy)=\infty$, 
 the local time may be a rough path of roughness $p>3$.

Having established the path integration of local time and the
corresponding convergence results, as a simple application, we can
easily prove a useful extension of It${\hat {\rm o}}$'s formula for
the L\'evy process when the function is less smooth:
 if $f:R\to
 R$ is an absolutely continuous function and has left
 derivative $\nabla^-f(x)$ being left continuous  and
  of bounded $q$-variation, where $1\leq
 q<3$, then  P-a.s.
 \begin{eqnarray}\label{may}
f(X_t)&=&f(X_0)+\int_0^t \nabla
^-f(X_s)dX_s-\int_{-\infty}^{\infty}\nabla ^-f(x)d_x
L_t^x\nonumber\\
&&+\sum_{0\leq s\leq t}[f(X_s)-f(X_{s-})-\Delta X_s\nabla
^-f(X_{s-})], \ \ 0\le t<\infty.
\end{eqnarray}
Here the path integral $\int_{-\infty}^{\infty}\nabla ^-f(x)d_x L_t^x$ is
a Lebesgue-Stieltjes integral when $q=1$, a Young integral when
$1<q<2$ and a Lyons' rough path integral when $2\leq q<3$
respectively. Needless to say that Tanaka's formula (\cite{tan}) and Meyer's formula
(\cite{meyer}, \cite{wang}) are very special cases of our formula when $q=1$. The investigation of
 It${\hat {\rm o}}$'s formula to less smooth functions is crucial and useful in many problems e.g. studying partial differential equations
 with singularities, the free boundary problem in American options,
 and certain stochastic differential equations.
 Time dependent cases for a continuous semimartingale $X_t$
 were investigated recently by Elworthy, Truman and Zhao (\cite{Zhao1}),
 Feng and Zhao (\cite{Zhao33}), where two-parameter Lebesgue-Stieltjes integrals and
 two-parameter Young integral
 were used respectively.  We would like to point out that a two-parameter rough path integration theory, which is important to the study of local times, and some other problems such as SPDEs, still remains open.

A part of the results about the rough path integral of local time
for a continuous semimartingale was announced in  Feng and Zhao
\cite{Feng3}. In summary, we have obtained complete results of the
variation and rough path of roughness $p$ for any $p\in (2,3)$ of local times for
any continuous semi-martingales and a class of L\'evy processes satisfying
(\ref{levynew36}) and $\sigma\neq 0$. There is no need to include a full
proof of the rough path result for
continuous semimartingales in this paper. We believe the reader can see easily that the proof is essentially
included in this paper, noticing the idea of
decomposing the local time to continuous and  discontinuous parts in \cite{Zhao33} and the key estimate (8) in
\cite{Feng3}.

\section {The $p$-variation of the local time of L\'evy processes}

Let $X_t$ be a one dimensional time homogeneous L\'evy process, and
$({\cal F} _t)_{t\geq 0}$ be generated by the sample paths $X_t$,
$p(\cdot)$ be a stationary $({\cal F} _t)$-Poisson process on
$R\setminus \{0\}$. From the well-known L\'evy-It$\hat{\rm o}$
decomposition theorem, we can write $X_t$ as follows:
\begin{eqnarray}\label{tanakaf}
X_t
%&=&X_0+ \sigma B_t+bt+\int_0^{t+}\int_{R\setminus \{0\}}y
%1_{\{|y|\geq 1\}} N_p(dsdy)+\int_0^{t+}\int_{R\setminus \{0\}}y
%1_{\{|y|< 1\}}
%\tilde N_p(dsdy)\\
 &:=&X_0+\sigma B_t+V_t+\tilde M_t,
\end{eqnarray}
where
\begin{eqnarray*}
V_t&:=&bt+\int_0^{t+}\int_{R\setminus \{0\}}y 1_{\{|y|\geq 1\}}
N_p(dsdy),\\
\tilde M_t&:=&\int_0^{t+}\int_{R\setminus \{0\}}y 1_{\{|y|< 1\}}
\tilde N_p(dsdy).
\end{eqnarray*}
Here, $N_p$ is the Poisson random measure of $p$, the compensator of
$p$ is of the form $\hat N(dsdy)=dsn(dy)$, where $n(dy)$ is the
L\'evy measure of process $X$. The compensated random measure
$\tilde N_p(t, U)=N_p(t, U)-\hat N_p(t, U)$ is an $({\cal F}
_t)$-martingale. Before we give our main theorem of this section, we
first give a $p$-moment estimate formula. This will be used in later
proofs.
\begin{lem}
%Assume $f$ is $({\cal F} _t)$-predictable, and satisfies
%$$E\bigg(\int_0^{t}\int_R |f(s,y,\omega)|^i
%n(dy)ds\bigg)<\infty, \ \ i=1,2,\cdots 4m.$$ Then
We have the p-moment estimate formula: for any $p\geq 1$,
\begin{eqnarray}
&&E\bigg(\sum_{0\leq s\leq t}|f(s,p_s(\omega),\omega)|\bigg)^p\nonumber\\
& \leq& 
%pE\bigg(\int_0^{t}\int_R
%|f(s,y,\omega)|n(dy)ds\bigg)^p
c_p\sum_{k=0}^mE\bigg(\int_0^{t}\int_R
|f(s,y,\omega)|^{2^k} n(dy)ds\bigg)^{p\over
{2^k}}
%&&+c_pE\bigg(\int_0^{t}\int_R |f(s,y,\omega)|^2
%n(dy)ds\bigg)^{p\over 2}+\cdots+c_pE\bigg(\int_0^{t}\int_R
%|f(s,y,\omega)|^{2m} n(dy)ds\bigg)^{p\over
%{2m}}\nonumber\\
+c_p\bigg(E\int_0^{t}\int_R |f(s,y,\omega)|^{2^{m+1}}
n(dy)ds\bigg)^{p\over {2^{m+1}}},
\end{eqnarray}
for a constant $c_p>0$. Here $m$ is the smallest integer such that
$2^{m+1}\geq p$.
\end{lem}
{\bf Proof}: From the definition of $N_p$, $\tilde N_p$, the
Burkholder-Davis-Gundy inequality and Jensen's inequality, we can
have the p-moment estimation:
\begin{eqnarray*}
 &&E\bigg(\sum_{0\leq s\leq
t}|f(s,p_s(\omega),\omega)|\bigg)^p\\
&=&E\bigg(\int_0^{t+}\int_R |f(s,y,\omega)|N_p(dsdy)\bigg)^p\\
&=&E\bigg(\int_0^{t}\int_R |f(s,y,\omega)|n(dy)ds+\int_0^{t+}\int_R |f(s,y,\omega)|\big(N_p(dsdy)-n(dy)ds\big)\bigg)^p\\
&\leq& pE\bigg(\int_0^{t}\int_R
|f(s,y,\omega)|n(dy)ds\bigg)^p+pE\bigg(\int_0^{t+}\int_R
|f(s,y,\omega)|\tilde
N_p(dsdy)\bigg)^p\\
&\leq&pE\bigg(\int_0^{t}\int_R |f(s,y,\omega)|n(dy)ds\bigg)^p+
c_pE\bigg(\int_0^{t}\int_R |f(s,y,\omega)|^2 n(dy)ds\bigg)^{p\over
2}\nonumber\\
&&+c_pE\bigg(\int_0^{t}\int_R |f(s,y,\omega)|^2 \tilde
N_p(dsdy)\bigg)^{p\over 2}\nonumber\\
&\leq&pE\bigg(\int_0^{t}\int_R |f(s,y,\omega)|n(dy)ds\bigg)^p+
c_pE\bigg(\int_0^{t}\int_R |f(s,y,\omega)|^2 n(dy)ds\bigg)^{p\over
2}\nonumber\\
&&+\cdots+c_pE\bigg(\int_0^{t}\int_R |f(s,y,\omega)|^{2^m}
n(dy)ds\bigg)^{p\over {2^m}}+c_p\left(E\bigg(\int_0^{t}\int_R
|f(s,y,\omega)|^{2^m} \tilde N_p(dsdy)\bigg)^2\right)^{{p\over
{2^{m+1}}}}\nonumber\\
&=&pE\bigg(\int_0^{t}\int_R |f(s,y,\omega)|n(dy)ds\bigg)^p+
c_pE\bigg(\int_0^{t}\int_R |f(s,y,\omega)|^2 n(dy)ds\bigg)^{p\over
2}\nonumber\\
&&+\cdots+c_pE\bigg(\int_0^{t}\int_R |f(s,y,\omega)|^{2^m}
n(dy)ds\bigg)^{p\over {2^m}}+c_p\bigg(E\int_0^{t}\int_R
|f(s,y,\omega)|^{{{2^{m+1}}}} n(dy)ds\bigg)^{{{p\over {2^{m+1}}}}},
%&=& pE\bigg(\int_0^{t}\int_R
%|f(s,y,\omega)|n(dy)ds\bigg)^p+c_p\sum_{k=1}^mE\bigg(\int_0^{t}\int_R
%|f(s,y,\omega)|^{2^k} n(dy)ds\bigg)^{p\over
%{2^k}}\nonumber\\
%&&+c_pE\bigg(\int_0^{t}\int_R |f(s,y,\omega)|^2
%n(dy)ds\bigg)^{p\over 2}+\cdots+c_pE\bigg(\int_0^{t}\int_R
%|f(s,y,\omega)|^{2m} n(dy)ds\bigg)^{p\over
%{2m}}\nonumber\\
%&&+c_p\bigg(E\int_0^{t}\int_R |f(s,y,\omega)|^{2^{m+1}}
%n(dy)ds\bigg)^{p\over {2^{m+1}}},
\end{eqnarray*}
where $m$ is the smallest integer such that $2^{m+1}\geq p$.
$\hfill\diamond$
\bigskip

 Recall the Tanaka formula for the L\'evy process
$X_t$ (c.f. \cite{apple}), we have
\begin{eqnarray}\label{local418}
L_t^a&=&(X_t-a)^+-(X_0-a)^+-\int_0^t 1_{\{X_{s-}>
a\}}dX_s\nonumber\\
&&+\sum_{0\leq s\leq t}[(X_{s-}-a)^+-(X_s-a)^++1_{\{X_{s-}>
a\}}\Delta X_s].
%&=&(X_t-a)^+-(X_0-a)^+-\int_0^t 1_{\{X_{s-}\geq a\}}d(cs+bB_s+\tilde
%M_s)\\
%&&+\sum_{0\leq s\leq t}[(X_{s-}-a)^+-(X_s-a)^+]1_{\{|\Delta
%X_s|\geq 1\}}\\
%&&+\sum_{0\leq s\leq t}[(X_{s-}-a)^+-(X_s-a)^++1_{\{X_{s-}\geq
%a\}}\Delta X_s]1_{\{|\Delta
%X_s|< 1\}}\\
\end{eqnarray}
Since
\begin{eqnarray*}
\sum_{0\leq s\leq t}1_{\{X_{s-}> a\}}\Delta X_s 1_{\{|\Delta
X_s|\geq 1\}}=\int_0^{t+} \int_{R\setminus \{0\}}1_{\{X_{s-}>a\}}y
1_{\{|y|\geq 1\}}N_p(dsdy),
\end{eqnarray*}
the following alternative form is often convenient
\begin{eqnarray}\label{levy3}
L_t^{a}=\varphi_t(a)-bI_t(a)-\sigma\hat
B^a_t+K_1(t,a)+K_2(t,a)+K_3(t,a),
\end{eqnarray}
where
\begin{eqnarray*}
\varphi_t(a)&:=&(X_t-a)^+-(X_0-a)^+,\ \ I_t(a):=\int_0^t
1_{\{X_{s-}> a\}}ds,\ \
\hat B^a_t:=\int_0^t 1_{\{X_{s-}>a\}} dB_s,\\
K_1(t,a)&:=&\sum_{0\leq s\leq t}[(X_{s-}-a)^+-(X_s-a)^+]1_{\{|\Delta
X_s|\geq 1\}},\\
K_2(t,a)&:=&\int_0^t\int_{R\setminus \{0\}}[(X_{s}-y-a)^+-(X_s-a)^+]1_{\{|y|< 1\}}\tilde N_p(dyds),\\
K_3(t,a)&:=&\int_0^t\int_R[(X_{s}-y-a)^+-(X_s-a)^++1_{\{X_{s}-y>
a\}}y]1_{\{|y|< 1\}}n(dy)ds.
\end{eqnarray*}
 For the convenience in what
follows in later part, we denote
\begin{eqnarray*}
J_1(s,a):=(X_{s-}-a)^+-(X_s-a)^+,\ \ J_2(s,a):=1_{\{X_{s-}>
a\}}\Delta X_s.
\end{eqnarray*}
Note we have the following important decompositions that we will use
often: for any $a_i<a_{i+1}$,
\begin{eqnarray}\label{levy57}
&&J_*(X_s,X_{s-},a_i,a_{i+1})\nonumber\\
&:=&J_1(s,a_{i+1})-J_1(s,a_i)\nonumber\\
&=&-(X_{s-}-a_i)1_{\{X_s\leq a_i\}}1_{\{a_i<
X_{s-}\leq a_{i+1}\}}-(a_{i+1}-a_i)1_{\{X_s\leq a_i\}}1_{\{X_{s-}>a_{i+1}\}}\nonumber\\
&&+(X_s-a_i)1_{\{a_i<X_{s}\leq a_{i+1}\}}1_{\{X_{s-}\leq a_i\}}
-(a_{i+1}-X_s)1_{\{a_i<X_{s}\leq a_{i+1}\}}1_{\{X_{s-}>
a_{i+1}\}}\nonumber\\
&&+(a_{i+1}-a_i)1_{\{X_{s}>a_{i+1}\}}1_{\{X_{s-}\leq a_i\}}
+(a_{i+1}-X_{s-})1_{\{X_{s}>a_{i+1}\}}1_{\{a_i< X_{s-}\leq
a_{i+1}\}}\nonumber\\
&&+ (X_s-X_{s-})1_{\{a_i<X_{s}\leq a_{i+1}\}}1_{\{a_i<X_{s-}\leq
a_{i+1}\}},
\end{eqnarray}
and
\begin{eqnarray}\label{levy9}
&&J^*(X_s,X_{s-},a_i,a_{i+1})\nonumber\\
&:=&[J_1(s,a_{i+1})-J_1(s,a_i)]+[J_2(s,a_{i+1})-J_2(s,a_i)]\nonumber\\
&=&-(X_{s}-a_i)1_{\{X_s\leq a_i\}}1_{\{a_i<
X_{s-}\leq a_{i+1}\}}-(a_{i+1}-a_i)1_{\{X_s\leq a_i\}}1_{\{X_{s-}>a_{i+1}\}}\nonumber\\
&&+(X_s-a_i)1_{\{a_i<X_{s}\leq a_{i+1}\}}1_{\{X_{s-}\leq a_i\}}
-(a_{i+1}-X_s)1_{\{a_i<X_{s}\leq a_{i+1}\}}1_{\{X_{s-}>
a_{i+1}\}}\nonumber\\
&&+(a_{i+1}-a_i)1_{\{X_{s}>a_{i+1}\}}1_{\{X_{s-}\leq a_i\}}
+(a_{i+1}-X_{s})1_{\{X_{s}>a_{i+1}\}}1_{\{a_i< X_{s-}\leq
a_{i+1}\}}.
\end{eqnarray}

The following theorem on the $p$-variation of local time in the
spatial variable is the main result of this section. We need to
consider the variation of each term in (\ref{levy3}). As the terms
related to the continuous part of the L\'evy process were considered
in \cite{Zhao33}, so the main difficulty is from the jumps related
to L\'evy process, especially small jumps. The main idea to deal
with the jump part is to use the $p$-moment estimate formula for the
jump part and change it to integration of $X_s$. Then we can use
occupation times formula, Jensen's inequality and Fubini theorem to
obtain the $p$-moment of the increment from $a_i$ to $a_{i+1}$.
%The
%decomposition (\ref{levy9}) of
%$J_1(s,a_{i+1})-J_1(s,a_{i})+J_2(s,a_{i+1})-J_2(s,a_{i})$ plays an
%important role in our proof of the desired increment estimates.
With the increment estimate, the $p$-variation can then be proved by
the classical approach of L\'evy.
\begin{thm}\label{newlem}
If $\sigma\neq 0$ and the L\'evy measure $n(dy)$ satisfies
\begin{eqnarray}\label{levy110}
\int_{R\setminus \{0\}}(|y|^{3\over 2}\wedge1)n(dy)<\infty,
\end{eqnarray}
then the local time $L_t^a$ of time homogeneous L\'evy process $X_t$
given by (\ref{tanakaf1}) is of bounded $p$-variation in $a$ for any
$t\geq 0$, for any $p>2$, almost surely, i.e.
\begin{eqnarray*}
\sup_{D(-\infty,\infty)}\sum_{i}|L_t^{a_{i+1}}-L_t^{a_{i}}|^p<\infty\
\ a.s.,
\end{eqnarray*}
where the supremum is taken over all finite partition on $R$,
$D(-\infty,\infty):=\{-\infty<a_0<a_1<\cdots<a_n<\infty\}$.
\end{thm}
{\bf Proof}: Will estimate every term in (\ref{levy3}).  First note that the function
$\varphi_t(a):=(X_t-a)^+-(X_0-a)^+$ is Lipschitz continuous in $a$
with Lipschitz constant 2. This implies that for any $p>2$ and
$a_i<a_{i+1}$,
\begin{eqnarray}\label{hz10}
| \varphi_t(a_{i+1})- \varphi_t(a_i)|^p\leq 2^p (a_{i+1}-a_i)^p.
\end{eqnarray}
Secondly, for the second term, by the occupation times formula,
Jensen's inequality and Fubini theorem,
%\begin{eqnarray}\label{hz11}
%|I_t(a_{i+1})-I_t(a_i)|^p=\left|\int_0^t 1_{\{a_i<X_s\leq
%a_{i+1}\}}ds\right|^p\leq c\int_0^t 1_{\{a_i<X_s\leq a_{i+1}\}}ds.
%\end{eqnarray}
\begin{eqnarray}\label{levy13a}
E|I_t(a_{i+1})-I_t(a_i)|^p
&=&E\left(\int_0^t 1_{\{a_i<X_s\leq a_{i+1}\}}ds\right)^{p}\nonumber\\
&=& {1\over {\sigma^{2p}}}(a_{i+1}-a_{i})^{p}E({1\over
a_{i+1}-a_{i}} \int _{a_i}^{a_{i+1}}L_t^xdx)^{p}
\nonumber\\
%&\leq & {1\over {\sigma^{2p}}}(a_{i+1}-a_{i})^{p-1}E{1\over
%a_{i+1}-a_{i}} \int _{a_i}^{a_{i+1}}(L_t^x)^{p}dx
%\nonumber\\
&\leq & {1\over {\sigma^{2p}}}(a_{i+1}-a_{i})^{p}\sup_{x}
E(L_t^x)^{p}.
\end{eqnarray}
By (\ref{local418}) and noting that $\sum\limits_{0\leq s\leq
t}[(X_{s-}-a)^+-(X_s-a)^++1_{\{X_{s-}> a\}}\Delta X_s]$
 is a decreasing process in $t$, we have
\begin{eqnarray*}
L_t^a\leq(X_t-a)^+-(X_0-a)^+-\int_0^t 1_{\{X_{s-}> a\}}dX_s.
\end{eqnarray*}
 Now
using the Burkholder-Davis-Gundy inequality again, we have
\begin{eqnarray*}
E(L_t^a)^{p}&\leq&{p}\Big[ E|X_t-X_0|^{p}+E|\int_0^t \sigma
1_{\{X_{s-}> a\}}dB_s|^{p}+
E|\int_0^t 1_{\{X_{s-}> a\}}dV_s|^{p}+E|\int_0^t 1_{\{X_{s-}> a\}}d{\tilde M}_s|^{p}\Big]\\
%&\leq& cE[|X_t-X_0|^{p\over 2}+t^{p\over 2}+<B>_t^{p\over 4}+<\tilde
%M>_t^{p\over 4}+|K(t,a)|^{p\over 2}]
%\\
&\leq& c\sigma^{p}t^{p\over 2}+cE(\int_0^t
|dV_s|)^{p}+c_p\sum_{k=1}^{m+1}\left[\int_0^{t}\int_{R}|y|^{2^k}
1_{\{|y|< 1\}}
n(dy)ds\right]^{p\over {2^{k}}}\\
%&&+c_p\bigg(\int_0^{t}\int_R
%|y|^{2^{m+1}} n(dy)ds\bigg)^{p\over {2^{m+2}}}\\
%&&+\cdots+c_p\bigg(\int_0^{t}\int_R |y|^{2m} n(dy)ds\bigg)^{p\over
%{4m}}+c_p\bigg(\int_0^{t}\int_R
%|y|^{4m} n(dy)ds\bigg)^{p\over {8m}}\\
&\leq & c( p, b, \sigma, t),
%<c_1(K,p),
\end{eqnarray*}
where $m>0$ is the smallest integer such that $2^{m+1}\geq p$, $c(
p, b, \sigma, t)$ is a universal constant depending on $p, b,
\sigma$, and $t$. By Jensen's inequality, we also have
\begin{eqnarray}
E(L_t^a)^{p\over 2} \leq c( p, b, \sigma, t).
\end{eqnarray}
This inequality will be used later.
 So
\begin{eqnarray}\label{levy11a}
E|I_t(a_{i+1})-I_t(a_i)|^p\leq c( p, b, \sigma,
t)(a_{i+1}-a_{i})^{p}.
\end{eqnarray}
%Secondly, by H\"older's inequality, as $V$ is of bounded variation,
%so
%\begin{eqnarray}\label{hz11}
%&&|\widehat V_t^{a_{i+1}}-\widehat V_t^{a_i}|^p\nonumber\\
%&\leq&|\int_0^t 1_{\{a_i<X_{s-}\leq a_{i+1}\}}|dV_s||^p\nonumber\\
%&\leq&c\int_0^t 1_{\{a_i<X_{s-}\leq a_{i+1}\}}|dV_s|,
%\end{eqnarray}
%where $c$ is a generic constant.
%Thirdly,
%\begin{eqnarray*}
%&&\sum_i\left|\int_0^t \int_{R\setminus \{0\}}1_{\{a_i<X_{s-}\leq
%a_{i+1}\}}x 1_{\{|x|\geq
%1\}}N_p(dsdx)\right|\\
%&\leq& \int_0^t\int_{R\setminus \{0\}} 1_{\{-N<X_{s-}\leq N\}}|x|
%1_{\{|x|\geq
%1\}} N_p(dsdx)\\
%&=&\sum_{0\leq s\leq t}|\Delta X_s|1_{\{|\Delta X_s|\geq 1\}}\\
%&<&\infty.
%\end{eqnarray*}
Thirdly, for the term $\hat B_t^a$,
\begin{eqnarray}\label{levy4}
E|\hat B_t^{a_{i+1}}-\hat B_t^{a_i}|^p
%&=&E|\int_0^t 1_{\{a_i<X_{s-}\leq a_{i+1}\}}dB_s|^p\nonumber\\
&=&E|\int_0^t 1_{\{a_i<X_{s}\leq a_{i+1}\}}dB_s|^p\nonumber\\
&\leq &c_pE\left(\int_0^t 1_{\{a_i<X_s\leq
a_{i+1}\}}ds\right)^{p\over
2}\nonumber\\
%&=&{c_p\over {\sigma^p}}E(\int _{a_i}^{a_{i+1}}L_t^xdx)^{p\over 2}\nonumber\\
%&=& c(p,\sigma)(a_{i+1}-a_{i})^{p\over 2}E({1\over a_{i+1}-a_{i}}
%\int _{a_i}^{a_{i+1}}L_t^xdx)^{p\over 2}
%\nonumber\\
%&\leq & c(p,\sigma)(a_{i+1}-a_{i})^{p\over 2}E{1\over a_{i+1}-a_{i}}
%\int _{a_i}^{a_{i+1}}(L_t^x)^{p\over 2}dx
%\nonumber\\
%&\leq & c(p,\sigma)(a_{i+1}-a_{i})^{p\over 2}\sup_{x}
%E(L_t^x)^{p\over 2}.
%\end{eqnarray}
%Here we used Burkholder-Davis-Gundy inequality, the occupation times
%formula, Jensen inequality and Fubini theorem.
% Hence it follows from Lemma \ref {lem2} that
%\begin{eqnarray}\label{levy5}
%E|\int_0^t 1_{\{a_i<X_{s-}\leq a_{i+1}\}}dB_s|^p
&\leq& c(t,p,\sigma)(a_{i+1}-a_{i})^{p\over 2}.
\end{eqnarray}
The last estimate can be obtained similarly to (\ref{levy13a}).
 About $K_1(t,a)$, it is easy to see that
\begin{eqnarray*}
&&|K_1(t,a_{i+1})-K_1(t,a_i)|\\
&\leq&\sum_{0\leq s \leq
t}|(X_{s-}-a_{i+1})^+-(X_s-a_{i+1})^+-(X_{s-}-a_{i})^++(X_s-a_{i})^+|1_{\{|\Delta X_s|\geq 1\}}\\
&\leq&2(a_{i+1}-a_i)\sum_{0\leq s \leq t}1_{\{|\Delta X_s|\geq 1\}}.
\end{eqnarray*}
So
\begin{eqnarray}\label{levy117}
&&E|K_1(t,a_{i+1})-K_1(t,a_i)|^p\leq C(a_{i+1}-a_i)^p.
\end{eqnarray}
About $K_2(t, a)$, with the decomposition (\ref{levy57}), we will
estimate the sum of each term for jumps $|\Delta X_s|<1$. There are
seven such terms.
%Recall (\ref{levy9}) and denote
%\begin{eqnarray*}
%J^*(X_s,X_{s-},a_i,a_{i+1})&=&J^*(X_s,X_{s}-\Delta X_s,a_i,a_{i+1})\\
%&:=&\Big[J_1(s,a_{i+1})-J_1(s,a_i)+J_2(s,a_{i+1})-J_2(s,a_i)\Big]1_{\{|\Delta
%X_s|<1\}}.
%\end{eqnarray*}
%We still write $K(t,a)$ like
%\begin{eqnarray*}
%K(t,a)&=&\sum_{0\leq s\leq t}[(X_{s-}-a)^+-(X_s-a)^++1_{\{X_{s-}\geq
%a\}}\Delta X_s]\\
%&=&\sum_{0\leq s\leq t} (J_1(s,a)+J_2(s,a))1_{\{|\Delta X_s|\geq
%1\}}+\sum_{0\leq s\leq t} (J_1(s,a)+J_2(s,a))1_{\{|\Delta X_s|<1\}}
%\end{eqnarray*}
%For the first sum, because
%\begin{eqnarray*}
%&&\sum_i\sum_{0\leq s\leq t} |J_2(s,a_{i+1})-J_2(s,a_i)|1_{\{|\Delta
%X_s|\geq 1\}}\\
%&\leq&\sum_i\sum_{0\leq s\leq t} 1_{\{a_i<X_{s-}\leq
%a_{i+1}\}}|\Delta X_s|1_{\{|\Delta X_s|\geq 1\}}\\
%&\leq&\sum_{0\leq s\leq t} |\Delta X_s|1_{\{|\Delta X_s|\geq 1\}}\\
%&<&\infty,
%\end{eqnarray*}
%so let's estimate the $p$-moment of $\sum\limits_{0\leq s\leq
%t}|J_1(s,a_{i+1})-J_1(s,a_i)|$. We will estimate every term of
%(\ref{levy7}).

 For the first term in (\ref{levy57}), by the $p$-moment estimate formula and occupation times formula, we have 
\begin{eqnarray}\label{new21}
&&E\left(\int_0^t\int_{R\setminus \{0\}}|X_{s-}-a_i|1_{\{X_s\leq a_i\}}1_{\{a_i<X_{s-}\leq a_{i+1}\}}1_{\{|\Delta X_s|< 1\}}
\tilde N_p(dyds)\right)^p\nonumber\\
&\leq&c_p \sum_{k=1}^m E\left(\int_0^t\int_{X_s-a_{i+1}}^{X_s-a_i}
|X_{s}-y-a_i|^{2^k}1_{\{X_s\leq a_i\}}1_{\{|y|< 1\}}
n(dy)ds\right)^{p\over {2^k}}\nonumber\\
%&&+\cdots +c_p E\left(\int_0^t\int_{X_s-a_{i+1}}^{X_s-a_i}
%|X_{s}-y-a_i|^{2m}1_{\{X_s\leq a_i\}}1_{\{|y|< 1\}}
%n(dy)ds\right)^{p\over {2m}}\\
&&+c_p \left(E\int_0^t\int_{X_s-a_{i+1}}^{X_s-a_i}
|X_{s}-y-a_i|^{{2^{m+1}}}1_{\{X_s\leq a_i\}}1_{\{|y|< 1\}}
n(dy)ds\right)^{p\over {2^{m+1}}}\nonumber\\
 &=&{ c_p\over {\sigma ^p}}
\sum_{k=1}^m
E\left(\int_{-\infty}^{a_i}L_t^x\int_{x-a_{i+1}}^{x-a_i}
|x-y-a_i|^{2^k} 1_{\{|y|<1\}}n(dy)dx\right)^{p\over
{2^k}}\nonumber\\
&&+{ c_p\over {\sigma
^p}}\left(E\int_{-\infty}^{a_i}L_t^x\int_{x-a_{i+1}}^{x-a_i}
|x-y-a_i|^{{2^{m+1}}} 1_{\{|y|<1\}}n(dy)dx\right)^{p\over
{{2^{m+1}}}},
\end{eqnarray}
where $m$ is the smallest integer such that $2^{m+1}\geq p$. 
In the following we will often use the following type of method to estimate
integrals with respect to the L\'evy measure: let
\begin{eqnarray*}
Q&:=&\int_{a_i-a_{i+1}}^0\int_{y+a_{i}}^{a_{i}}
(x-y-a_i)dx1_{\{|y|< 1\}}n(dy)\\
&=&{1\over 2}\int_{a_i-a_{i+1}}^0|y|^21_{\{|y|< 1\}}n(dy).
\end{eqnarray*}
Then, $1_{\{|y|<1\}}{1\over Q}(x-y-a_i)dxn(dy)$ is a probability measure on
$\{(x,y): y+a_i\leq x\leq a_i, (-1)\wedge (a_i-a_{i+1})\leq y\leq
0\}$. So by Jensen's inequality, we have that
\begin{eqnarray}\label{1june2009}
&&E\left(\int_{a_i-a_{i+1}}^0\int_{y+a_{i}}^{a_{i}}L_t^x
(a_{i+1}-a_i)(x-y-a_i)dx1_{\{|y|< 1\}}n(dy)\right)^{p\over 2}\nonumber\\
&\leq & Q^{p\over 2}E\left({1\over
Q}\int_{a_i-a_{i+1}}^0\int_{y+a_{i}}^{a_{i}}(L_t^x)^{p\over 2}
(a_{i+1}-a_i)^{{p\over 2}}(x-y-a_i)dx1_{\{|y|<
1\}}n(dy)\right)\nonumber\\
&\leq&c(\sigma,t,p)|a_{i+1}-a_i|^{{p\over 2}}\sup_x E(L_t^x)^{p\over
2}\left (\int_{a_i-a_{i+1}}^0|y|^2 1_{\{|y|< 1\}}n(dy)\right )^{p\over 2}.
\end{eqnarray}
Similarly one can estimate
\begin{eqnarray}
&&
E\left (\int_{-\infty}^{a_i-a_{i+1}}\int_{y+a_{i}}^{y+a_{i+1}}L_t^xdx
|y|^21_{\{|y|<
1\}}n(dy)\right )^{p\over 2}\nonumber
\\
&\leq & c(\sigma ,t,p) |a_{i+1}-a_i|^{p\over 2}\sup_x E(L_t^x)^{p\over
2}\left (\int_{-\infty}^{a_i-a_{i+1}}|y|^2 1_{\{|y|<
1\}}n(dy)\right )^{p\over 2}.
\end{eqnarray}
Then we can estimate each term in (\ref{new21}). When $k=1$, we
change the orders of the integration and use Jensen's inequality to have 
\begin{eqnarray}\label{new19}
&&E\left(\int_{-\infty}^{a_i}L_t^x\int_{x-a_{i+1}}^{x-a_i}
|x-y-a_i|^{2} 1_{\{|y|<1\}}n(dy)dx\right)^{p\over {2}}\nonumber\\
 &=&E\Bigg(\int_{-\infty}^{a_i-a_{i+1}}\int_{y+a_{i}}^{y+a_{i+1}}L_t^x
(x-y-a_i)^21_{\{|y|<
1\}}dxn(dy)\nonumber\\
&&\hskip 1.5cm+\int_{a_i-a_{i+1}}^0\int_{y+a_{i}}^{a_{i}}L_t^x
(x-y-a_i)(x-y-a_i)1_{\{|y|< 1\}}dxn(dy)\Bigg)^{p\over 2}\nonumber\\
&\leq &
E\Bigg(\int_{-\infty}^{a_i-a_{i+1}}\int_{y+a_{i}}^{y+a_{i+1}}L_t^xdx
|y|^21_{\{|y|<
1\}}n(dy)\nonumber\\
&&\hskip 1.5cm+\int_{a_i-a_{i+1}}^0\int_{y+a_{i}}^{a_{i}}L_t^x
(a_{i+1}-a_i)(x-y-a_i)dx1_{\{|y|< 1\}}n(dy)\Bigg)^{p\over 2}\\
&\leq & c(\sigma, t,p) |a_{i+1}-a_i|^{p\over 2}\sup_x E(L_t^x)^{p\over
2}\left (\int_{-\infty}^{a_i-a_{i+1}}|y|^2 1_{\{|y|<
1\}}n(dy)\right )^{p\over 2}\nonumber\\
&&+c(\sigma,t,p)|a_{i+1}-a_i|^{{p\over 2}}\sup_x E(L_t^x)^{p\over
2}\left (\int_{a_i-a_{i+1}}^0|y|^2 1_{\{|y|< 1\}}n(dy)\right )^{p\over 2}\nonumber\\
 &\leq
&c(\sigma, t,p)|a_{i+1}-a_i|^{p\over 2}\nonumber.
\end{eqnarray}

For the term when $2\leq k\leq m$,
\begin{eqnarray*}
&&E\left(\int_{-\infty}^{a_i}L_t^x\int_{x-a_{i+1}}^{x-a_i}
|x-y-a_i|^{2^k} 1_{\{|y|<1\}}n(dy)dx\right)^{p\over {2^k}}\\
&\leq &E\left(\int_{-\infty}^{a_i}L_t^x\int_{x-a_{i+1}}^{x-a_i}
|x-y-a_i|^{2^{k-1}}|y|^{2^{k-1}} 1_{\{|y|<1\}}n(dy)dx\right)^{p\over {2^k}}\\
&\leq &(a_{i+1}-a_i)^{{1\over
2}p}E\left(\int_{-\infty}^{a_i}L_t^x\int_{x-a_{i+1}}^{x-a_i}
|y|^{2^{k-1}} 1_{\{|y|<1\}}n(dy)dx\right)^{p\over {2^k}}\\
&\leq &c(t,p)|a_{i+1}-a_i|^{p\over 2}.
\end{eqnarray*}
Actually we can see that
$\int_{-\infty}^{a_i}\int_{x-a_{i+1}}^{x-a_i} |y|^{2^{k-1}}
1_{\{|y|<1\}}n(dy)dx<\infty$ by using the same method as in
(\ref{new19}).
 For the last term in (\ref {new21}), similarly, we have
\begin{eqnarray*}
&&\left(E\int_{-\infty}^{a_i}L_t^x\int_{x-a_{i+1}}^{x-a_i}
|x-y-a_i|^{2^{m+1}} 1_{\{|y|<1\}}n(dy)dx\right)^{p\over {2^{m+1}}}\\
&\leq&\left(E\int_{-\infty}^{a_i}L_t^x\int_{x-a_{i+1}}^{x-a_i}
|a_{i+1}-a_i|^{2^m}|y|^{2^m} 1_{\{|y|<1\}}n(dy)dx\right)^{p\over
{2^{m+1}}}\\
&\leq&|a_{i+1}-a_i|^{p\over
2}\left(E\int_{-\infty}^{a_i}L_t^x\int_{x-a_{i+1}}^{x-a_i} |y|^{2^m}
1_{\{|y|<1\}}n(dy)dx\right)^{p\over
{2^{m+1}}}\\
&\leq &c(t,p)|a_{i+1}-a_i|^{p\over 2}.
\end{eqnarray*}
We can see that the key point is to estimate the term when $k=1$
because the higher order term can always be dealt by the above
method. We can use the similar method to deal with other terms and
derive that
\begin{eqnarray}\label{levy10}
&&E|K_2(t,a_{i+1})-K_2(t,a_i)|^p\leq
c(t,\sigma,p)|a_{i+1}-a_i|^{p\over 2}.
\end{eqnarray}
About $K_3(t,a)$, with the decomposition (\ref{levy9}), we will
estimate the sum of each term for jumps $|\Delta X_s|<1$. There are
six such terms.

 For the first term in (\ref{levy9}), by the
$p$-moment estimate formula and occupation times formula, we have
\begin{eqnarray*}
&&E\left(\int_0^t\int_{-\infty}^\infty |X_{s}-a_i|1_{\{X_s\leq a_i\}}1_{\{a_i<X_{s}-y\leq a_{i+1}\}}1_{\{|y|< 1\}}n(dy)ds\right)^p\\
 &=&E\left(\int_0^t\int_{X_s-a_{i+1}}^{X_s-a_i}
|X_{s}-a_i|1_{\{X_s\leq a_i\}}1_{\{|y|< 1\}}
n(dy)ds\right)^{p}\\
&=&{ 1\over {\sigma ^{2p}}}
E\left(\int_{-\infty}^{a_i}L_t^x\int_{x-a_{i+1}}^{x-a_i} |x-a_i|
1_{\{|y|< 1\}}n(dy)dx\right)^{p}.
\end{eqnarray*}
Now we change the orders of the integration and use Jensen's
inequality, we have
\begin{eqnarray}\label{levynew24}
&&E\left(\int_0^t\int_{-\infty}^\infty |X_{s}-a_i|1_{\{X_s\leq a_i\}}1_{\{a_i<X_{s}-y\leq a_{i+1}\}}1_{\{|y|< 1\}}n(dy)ds\right)^p\nonumber\\
 &\leq &{1\over {\sigma^{2p}}}
E\Bigg(\int_{-\infty}^{a_i-a_{i+1}}\int_{y+a_{i}}^{y+a_{i+1}}L_t^x
(a_i-x)1_{\{|y|<
1\}}dxn(dy)\nonumber\\
&&\hskip 1.5cm+\int_{a_i-a_{i+1}}^0\int_{y+a_{i}}^{a_{i}}L_t^x
(a_i-x)^{1\over 2}(a_i-x)^{1\over 2} 1_{\{|y|< 1\}}dxn(dy)\Bigg)^{p}\\
&\leq &{1\over {\sigma^{2p}}}
E\Bigg(\int_{-\infty}^{a_i-a_{i+1}}\int_{y+a_{i}}^{y+a_{i+1}}L_t^xdx
|y|1_{\{|y|<
1\}}n(dy)\nonumber\\
&&\hskip 1.5cm+\int_{a_i-a_{i+1}}^0\int_{y+a_{i}}^{a_{i}}L_t^x
(a_{i+1}-a_i)^{1\over 2}(a_i-x)^{1\over 2}dx1_{\{|y|< 1\}}n(dy)\Bigg)^{p}\nonumber\\
&\leq & c(\sigma) |a_{i+1}-a_i|^{{1\over 2}p}\sup_x
E(L_t^x)^{p}\left (\int_{-\infty}^{a_i-a_{i+1}}|y|^{{1\over 2}+1} 1_{\{|y|<
1\}}n(dy)\right )^p\nonumber\\
&&+c(\sigma)|a_{i+1}-a_i|^{{1\over 2}p}\sup_x E(L_t^x)^{p}\left (\int_{a_i-a_{i+1}}^0|y|^{3\over 2} 1_{\{|y|< 1\}}n(dy)\right )^p\nonumber\\
&\leq &c(t,\sigma,p)|a_{i+1}-a_i|^{{1\over 2}p}\nonumber.
\end{eqnarray}
We can use the similar method to deal with other terms. In the following, we will only sketch the estimate without giving great details.  \\
2) For the second term, we have
\begin{eqnarray*}
&&E\left(\int_0^t\int_{-\infty}^{\infty}(a_{i+1}-a_i)1_{\{X_s\leq
a_i\}}1_{\{X_{s}-y>
a_{i+1}\}}1_{\{|y|< 1\}}n(dy)ds\right)^{p}\\
&=&{1\over
{\sigma^{2p}}}(a_{i+1}-a_i)^{p}E\left(\int_{-\infty}^{a_i}L_t^x\int_{-\infty}^{x-
a_{i+1}}1_{\{|y|< 1\}}n(dy)dx\right)^{p}\\
&=&c(\sigma)(a_{i+1}-a_i)^{p}E\left(\int_{-\infty}^{a_i-a_{i+1}}\left(\int_{y+a_{i+1}}^{a_i}L_t^xdx\right)1_{\{|y|< 1\}}n(dy)\right)^{p}\\
&\leq &c(\sigma)|a_{i+1}-a_i|^{{1\over 2}p}\sup_x E(L_t^x)^{p}\left (\int_{-\infty}^{a_i-a_{i+1}}|y|^{{1\over 2}+1} 1_{\{|y|< 1\}}n(dy)\right )^p\\
&\leq &c(t,\sigma,p)|a_{i+1}-a_i|^{{1\over 2}p}.
\end{eqnarray*}
%We need
%\begin{eqnarray*}
% M_2:=\int_{-\infty}^{a}U(t,x)\Big(n({-\infty},(x-
%a-\varepsilon)\Big)^{p}dx\leq C(a)\varepsilon^{-{p\over 2}}
%\end{eqnarray*}

3) For the third term, we have
\begin{eqnarray}\label{levynew25}
&& E\left(\int_0^t \int_{-\infty}^\infty|X_s-a_i| 1_{\{a_i<
X_{s}\leq a_{i+1}\}}1_{\{X_{s}-y\leq a_i\}}1_{\{|y|< 1\}}n(dy)ds\right)^{p}\nonumber\\
&=& {1\over {\sigma^{2p}}} E\left(\int_{a_i}^{a_{i+1}}L_t^x\int_{x-a_i}^\infty |x-a_i|1_{\{|y|< 1\}}n(dy)dx\right)^{p}\nonumber\\
&\leq & {1\over {\sigma^{2p}}}
E\Bigg(\int_{0}^{a_{i+1}-a_i}\int_{a_i}^{a_i+y}L_t^x
(a_{i+1}-a_i)^{1\over 2} (x-a_i)^{1\over 2} dx1_{\{|y|< 1\}}n(dy)\nonumber\\
&&+\int_{a_{i+1}-a_i}^\infty\int_{a_i}^{a_{i+1}}L_t^x dx|y|1_{\{|y|< 1\}}n(dy)\Bigg)^{p}\\
&\leq &c(\sigma,p)|a_{i+1}-a_i|^{{1\over 2}p}\sup_x
E(L_t^x)^{p}\left (\int_{0}^{\infty}|y|^{3\over 2}1_{\{|y|< 1\}}n(dy)\right )^p\nonumber\\
&\leq &c(t,\sigma,p)|a_{i+1}-a_i|^{{1\over 2}p}\nonumber.
\end{eqnarray}
4) The fourth term is symmetric to the third term, so by a similar
computation, we have
\begin{eqnarray*}
&&E\left(\int_0^t\int_{-\infty}^\infty|a_{i+1}-X_s|1_{\{a_i<
X_{s}\leq a_{i+1}\}}1_{\{X_{s}-y>a_{i+1}\}}1_{\{|y|<
1\}}n(dy)ds\right)^p
%\\
%&\leq & c_p E\left(\int_0^t\int_{-\infty}^\infty
%|a_{i+1}-X_s|^21_{\{a_i<X_{s}\leq a_{i+1}\}}1_{\{X_{s}-y>
%a_{i+1}\}}1_{\{|y|< 1\}}n(dy)ds\right)^{p\over 2}\\
%&&+ p E\left(\int_0^t\int_{-\infty}^\infty |a_{i+1}-X_s|1_{\{a_i\leq
%X_{s}\leq a_{i+1}\}}1_{\{X_{s}-y>
%a_{i+1}\}}1_{\{|y|< 1\}}n(dy)ds\right)^p\\
%&=& {c_p\over {\sigma^{p}}}
%E\left(\int_{a_i}^{a_{i+1}}L_t^x\int_{-\infty}^{x-a_{i+1}}
%|a_{i+1}-x|^21_{\{|y|< 1\}}n(dy)dx\right)^{p\over 2}\\
%&&+{p\over {\sigma^{2p}}}
%E\left(\int_{a_i}^{a_{i+1}}L_t^x\int_{-\infty}^{x-a_{i+1}}
%|a_{i+1}-x|1_{\{|y|< 1\}}n(dy)dx\right)^{p}\\
%&\leq &c(\sigma,p) E\Bigg(\int_{a_i-a_{i+1}}^0
%\int_{y+a_{i+1}}^{a_{i+1}}L_t^x
%|a_{i+1}-a_i|^2dx1_{\{|y|< 1\}}n(dy)\\
%&&+\int_{-\infty}^{a_i-a_{i+1}} \int^{a_{i+1}}_{a_{i}}L_t^x
%|a_{i+1}-a_i||y|dx1_{\{|y|< 1\}}n(dy)\Bigg)^{p\over 2}\\
%&&+c(\sigma,p) E\Bigg(\int_{a_i-a_{i+1}}^0
%\int_{y+a_{i+1}}^{a_{i+1}}L_t^x
%|a_{i+1}-a_i|dx1_{\{|y|< 1\}}n(dy)\\
%&&+\int_{-\infty}^{a_i-a_{i+1}} \int^{a_{i+1}}_{a_{i}}L_t^x|y|dx1_{\{|y|< 1\}}n(dy)\Bigg)^{p}\\
%&\leq &c(\sigma,p)|a_{i+1}-a_i|^{p}\Big(\sup_x E(L_t^x)^{p\over
%2}+\sup_x
%E(L_t^x)^{p}\Big)\int_{-\infty}^0|y|1_{\{|y|< 1\}}n(dy)\\
\leq c(t,\sigma,p)|a_{i+1}-a_i|^{{1\over 2}p}.
%&\leq & c(t,b,p) \int_0^t E\left(\int_{-\infty}^{X_s-a_{i+1}}
%|a_{i+1}-X_s|^21_{\{a_i<X_{s}\leq a_{i+1}\}}n(dy)\right)^{p\over 2}ds\\
%&&+ c(t,b,p) \int_0^t E\left(\int_{-\infty}^{X_s-a_{i+1}}
%|a_{i+1}-X_s|1_{\{a_i<X_{s}\leq a_{i+1}\}}n(dy)\right)^{p}ds\\
%&=&c(t,b,p)\int_{a_i}^{a_{i+1}}U(t,x)|a_{i+1}-x|^p\left(\int_{-\infty}^{x-a_{i+1}}
%n(dy)\right)^{p}dx\\
%&&+c(t,b,p)\int_{a_i}^{a_{i+1}}U(t,x)|a_{i+1}-x|^p\left(\int_{-\infty}^{x-a_{i+1}}
%n(dy)\right)^{p\over 2}dx\\
%&=&c(t,b,p)\int_{a_i}^{a_{i+1}}U(t,x)|a_{i+1}-x|^p\Big(n({-\infty},{x-a_{i+1}})\Big)^{p}dx\\
%&&+c(t,b,p)\int_{a_i}^{a_{i+1}}U(t,x)|a_{i+1}-x|^p\Big(n({-\infty},{x-a_{i+1}})\Big)^{p\over 2}dx\\
%&\leq
%&c(t,b,p)(a_{i+1}-a_i)^p\left[\int_{a_i}^{a_{i+1}}U(t,x)\left(\int_{-\infty}^{x-a_{i+1}}
%n(dy)\right)^{p}dx+\int_{a_i}^{a_{i+1}}U(t,x)\left(\int_{-\infty}^{x-a_{i+1}}
%n(dy)\right)^{p\over 2}dx\right].
\end{eqnarray*}
%We need
%\begin{eqnarray*}
%M_4:=\int_{a}^{a+\varepsilon}U(t,x)|a+\varepsilon-x|^p\Big(n({-\infty},{x-a-\varepsilon})\Big)^{p}dx\leq
%\varepsilon^{{p\over 2}}.
%\end{eqnarray*}
%\begin{eqnarray*}
%M_5&:=&\int_{a_i}^{a_{i+1}}U(t,x)\left(\int_{-\infty}^{x-a_{i+1}}
%n(dy)\right)^{p}dx+\int_{a_i}^{a_{i+1}}U(t,x)\left(\int_{-\infty}^{x-a_{i+1}}
%n(dy)\right)^{p\over 2}dx\\
%&<&\infty.
%\end{eqnarray*}
5) For the fifth term, as it is symmetric to the second term, so we
can use a similar computation to have
\begin{eqnarray*}
E\left(\int_0^t\int_{-\infty}^\infty(a_{i+1}-a_i)1_{\{X_{s}>
a_{i+1}\}}1_{\{X_{s}-y\leq a_i\}}1_{\{|y|< 1\}}n(dy)ds\right)^p
%\\
%&\leq& (a_{i+1}-a_i)^p \Bigg[c_p E
%\left(\int_0^t\int_{-\infty}^{\infty}1_{\{X_s>
%a_{i+1}\}}1_{\{X_{s}-y\leq
%a_{i}\}}1_{\{|y|< 1\}}n(dy)ds\right)^{p\over 2}\\
%&&\hskip 2.5cm+pE\left(\int_0^t\int_{-\infty}^{\infty}1_{\{X_s>
%a_{i+1}\}}1_{\{X_{s}-y\leq
%a_{i}\}}1_{\{|y|< 1\}}n(dy)ds\right)^{p}\Bigg]\\
%&=&(a_{i+1}-a_i)^p \Bigg[{c_p\over {\sigma^{p}}} E
%\left(\int_{a_{i+1}}^{\infty}L_t^x\int_{x-a_i}^\infty 1_{\{|y|<
%1\}}n(dy)dx\right)^{p\over 2}\\
%&&\hskip 2.5cm+ {p\over {\sigma^{2p}}}E
%\left(\int_{a_{i+1}}^{\infty}L_t^x\int_{x-a_i}^\infty
%1_{\{|y|< 1\}}n(dy)dx\right)^{p}\Bigg]\\
%&\leq &c(\sigma,p)(a_{i+1}-a_i)^p \Bigg[ E
%\left(\int_{a_{i+1}-a_i}^{\infty}\int_{a_{i+1}}^{y+a_i}L_t^xdx
%1_{\{|y|<
%1\}}n(dy)\right)^{p\over 2}\\
%&&\hskip 3.5cm+E
%\left(\int_{a_{i+1}-a_i}^{\infty}\int_{a_{i+1}}^{y+a_i}L_t^xdx
%1_{\{|y|< 1\}}n(dy)\right)^{p}\Bigg]\\
%&\leq &c(\sigma,p)|a_{i+1}-a_i|^{p}\Big(\sup_x E(L_t^x)^{p\over
%2}+\sup_x
%E(L_t^x)^{p}\Big)\int_{a_{i+1}-a_i}^{\infty}|y|1_{\{|y|< 1\}}n(dy)\\
\leq c(t,\sigma,p)|a_{i+1}-a_i|^{{1\over 2}p}.
%&\leq& c(t,b,p)(a_{i+1}-a_i)^p \Bigg[ \int_0^t E\left(\int_{X_{s}-
%a_{i}}^\infty 1_{\{X_s>a_{i+1}\}}n(dy)\right)^{p\over 2}ds\\
%&&+\int_0^t E\left(\int_{X_{s}-
%a_{i}}^\infty 1_{\{X_s>a_{i+1}\}}n(dy)\right)^{p}ds\Bigg]\\
%&=&c(t,b,p)(a_{i+1}-a_i)^p\Bigg[\int_{a_{i+1}}^{\infty}U(t,x)\Big(n(x-
%a_{i},\infty)\Big)^{p\over
%2}dx+\int_{a_{i+1}}^{\infty}U(t,x)\Big(n(x-
%a_{i},\infty)\Big)^{p}dx\Bigg].
\end{eqnarray*}
%We need
%\begin{eqnarray*}
%M_5:=\int_{x-(a+\varepsilon)}^{\infty}U(t,x)\Big(n(x-a,\infty)\Big)^{p}dx\leq
%\varepsilon^{-{p\over 2}}.
%\end{eqnarray*}
%\begin{eqnarray*}
%M_6&:=&\int_{a_{i+1}}^{\infty}U(t,x)\left(\int_{x- a_{i}}^\infty
%n(dy)\right)^{p\over
%2}dx+\int_{a_{i+1}}^{\infty}U(t,x)\left(\int_{x- a_{i}}^\infty
%n(dy)\right)^{p}dx\\
%&<&\infty.
%\end{eqnarray*}
6) The last term is symmetric to the first term, so by a similar
computation, we have
\begin{eqnarray*}
E\left(\int_0^t\int_{-\infty}^\infty |a_{i+1}-X_{s}|1_{\{X_{s}>
a_{i+1}\}}1_{\{a_i<X_{s}-y\leq a_{i+1}\}}1_{\{|y|< 1\}}\right)^p
%\int_{x-a_{i+1}}^{x-a_i} |x-a_{i+1}|1_{\{|y|< 1\}}n(dy)dx\right)^{p}\Bigg]\\
%&\leq &c(\sigma,p) E\Bigg(\int_0^{a_{i+1}-a_i}
%\int_{a_{i+1}}^{a_{i+1}+y}L_t^x |a_{i+1}-a_i|^2dx1_{\{|y|<
%1\}}n(dy)\\
%&&\hskip 2cm+\int_{a_{i+1}-a_i}^\infty
%\int^{y+a_{i+1}}_{y+a_{i}}L_t^x |y|^2dx1_{\{|y|< 1\}}n(dy)\Bigg)^{p\over 2}\\
%&&+c(\sigma,p) E\Bigg(\int_0^{a_{i+1}-a_i}
%\int_{a_{i+1}}^{a_{i+1}+y}L_t^x |a_{i+1}-a_i|dx1_{\{|y|<
%1\}}n(dy)\\
%&&\hskip 2cm+\int_{a_{i+1}-a_i}^\infty
%\int^{y+a_{i+1}}_{y+a_{i}}L_t^x |y|dx1_{\{|y|< 1\}}n(dy)\Bigg)^{p}\\
%&\leq &c(\sigma,p)
%E\Bigg(\int_{-\infty}^{a_i-a_{i+1}}\int_{y+a_{i}}^{y+a_{i+1}}L_t^xdx
%|y|^21_{\{|y|<
%1\}}n(dy)\\
%&&\hskip 1.5cm+\int_{a_i-a_{i+1}}^0\int_{y+a_{i}}^{a_{i}}L_t^x
%|a_{i+1}-a_i|^2dx1_{\{|y|< 1\}}n(dy)\Bigg)^{p\over 2}\\
%&&+c(\sigma,p)
%E\Bigg(\int_{-\infty}^{a_i-a_{i+1}}\int_{y+a_{i}}^{y+a_{i+1}}L_t^xdx
%|y|1_{\{|y|<
%1\}}n(dy)\\
%&&\hskip 1.5cm+\int_0^{a_{i+1}-a_i}\int_{y+a_{i}}^{a_{i}}L_t^x
%|a_{i+1}-a_i|dx1_{\{|y|< 1\}}n(dy)\Bigg)^{p}\\
%&\leq &c(\sigma,p)|a_{i+1}-a_i|^{p}\sup_x E(L_t^x)^{p\over
%2}\int_0^{a_{i+1}-a_i}|y| 1_{\{|y|< 1\}}n(dy)\\
%&&+c(\sigma,p) |a_{i+1}-a_i|^{p\over 2}\sup_x E(L_t^x)^{p\over
%2}\int^{\infty}_{a_{i+1}-a_i}|y|^2 1_{\{|y|<
%1\}}n(dy)\\
%&&+ c(\sigma,p) |a_{i+1}-a_i|^{p}\sup_x
%E(L_t^x)^{p}\int^{\infty}_0|y| 1_{\{|y|< 1\}}n(dy)\\
\leq  c(t,\sigma,p)|a_{i+1}-a_i|^{{1\over 2}p}.
\end{eqnarray*}
So we have
\begin{eqnarray}\label{levy59}
E|K_3(t,a_{i+1})-K_3(t,a_i)|^p\leq
c(t,\sigma,p)|a_{i+1}-a_i|^{{1\over 2}p}.
\end{eqnarray}

 Now we use Proposition 4.1.1 in \cite{terry}
($i=1, \gamma>p-1$), for any finite partition $\{a_l\}$ of $[a,b]$
\begin{eqnarray*}
\sup_D\sum_l |\hat B_t^{a_{l+1}}-\hat B_t^{a_l}|^p\leq
c(p,\gamma)\sum_{n=1}^{\infty}n^{\gamma}\sum_{k=1}^{2^n}|\hat
B_t^{a_k^n}-\hat B_t^{a_{k-1}^n}|^p,
\end{eqnarray*}
where
\begin{eqnarray*}
a_k^n=a+{k\over {2^n}}(b-a),\ k=0,1,\cdots,2^n.
\end{eqnarray*}
The key point here is that the right hand side does not depend on
the partition D. We take the expectation and use (\ref{levy4}), it
follows that
\begin{eqnarray*}
E\sum_{n=1}^{\infty}n^{\gamma}\sum_{k=1}^{2^n}|\hat B_t^{a_k^n}-\hat B_t^{a_{k-1}^n}|^p&=&\sum_{n=1}^{\infty}n^{\gamma}\sum_{k=1}^{2^n}E|\hat B_t^{a_k^n}-\hat B_t^{a_{k-1}^n}|^p\\
&\leq& c\sum_{n=1}^{\infty}n^{\gamma}({{b-a}\over {2^n}})^{{p\over
2}-1}<\infty,
\end{eqnarray*}
as $p>2$. Therefore
\begin{eqnarray*}
\sum_{n=1}^{\infty}n^{\gamma}\sum_{k=1}^{2^n}|\hat B_t^{a_k^n}-\hat
B_t^{a_{k-1}^n}|^p<\infty\ a.s.
\end{eqnarray*}
It turns out that for any interval $[a,b]\subset R$
\begin{eqnarray}
\sup_D\sum_l |\hat B_t^{a_{l+1}}-\hat B_t^{a_l}|^p<\infty\ a.s.
\end{eqnarray}
But we know $L_t^{a}$ has a compact support $[-K,K]$ in $a$. So for
the partition $D:=D_{-K,K}=\{-K=a_0<a_1<\cdots<a_r=K\}$, we obtain
\begin{eqnarray}\label{levy17}
\sup_D\sum_l |\hat B_t^{a_{i+1}}-\hat B_t^{a_i}|^p<\infty\ a.s.
\end{eqnarray}
In the same way, from (\ref{levy11a}),
(\ref{levy117}),(\ref{levy10}), (\ref{levy59}) we can prove that
\begin{eqnarray}
&&\sup_{D}\sum_i |I_t^{a_{i+1}}-I_t^{a_i}|^p<\infty\
a.s.,\label{levy18}\\
&& \sup_{D}\sum_i |K_1(t,{a_{i+1}})-K_1(t,{a_i})|^p<\infty\
a.s.,\label{24a} \\
&& \sup_{D}\sum_i |K_2(t,{a_{i+1}})-K_2(t,{a_i})|^p<\infty\
a.s.\label{levy19} \\
&& \sup_{D}\sum_i |K_3(t,{a_{i+1}})-K_3(t,{a_i})|^p<\infty\
a.s.\label{levy20}
\end{eqnarray}
On the other hand, it is easy to see from (\ref{hz10}) that
\begin{eqnarray} \label{hz13}
&&\sum_i|\varphi_t(a_{i+1})-\varphi_t(a_i)|^p\leq
2^p\sum_i(a_{i+1}-a_i)^p \leq 2^p[\sum_i(a_{i+1}-a_i)]^p=
2^p(b-a)^p.
%&& \sum_i|I_t(a_{i+1})-I_t(a_i)|^p\leq c\int_0^t 1_{\{a<X_s\leq
%b\}}ds\leq ct.\label{levy24}
\end{eqnarray}
Then from (\ref{levy3}), (\ref{levy17}), (\ref{levy18}),(\ref{24a}),
(\ref{levy19}) (\ref{levy20}),  and (\ref{hz13}), it turns out that
\begin{eqnarray*}
\hskip 5cm \sup_{D}\sum_i |L_t^{a_{i+1}}-L_t^{a_i}|^p<\infty\ \ \ \
\ \ a.s.\hskip 5cm\hfill\diamond
\end{eqnarray*}
\begin{rmk}
From (\ref{levynew24}) and (\ref{levynew25}), we can see easily that
if we require the following slightly stronger condition on the
L\'evy measure
\begin{eqnarray}\label{levynew34}
\int_{R\setminus \{0\}}(|y|^{{3\over 2}-\xi}\wedge1)n(dy)<\infty,
\end{eqnarray}
for a $\xi\in (0,{1\over 2}]$, then for any $p\geq 1$,
\begin{eqnarray}\label{2june2009}
E|K_3(t,a_{i+1})-K_3(t,a_i)|^p\leq
c(t,\sigma,p)|a_{i+1}-a_i|^{({1\over 2}+\xi)p}.
\end{eqnarray}
This estimate will be used in the construction of the geometric
rough path where (\ref{levy59}) is not adequate.
\end{rmk}

%
% In the following, we assume
%\begin{eqnarray}\label{levy110}
%\int_{R\setminus \{0\}}(|y|^{4\over 3}\wedge1)n(dy)<\infty.
%\end{eqnarray}
%Under (\ref{levy110}), we will have more delicate estimate of the
%increment of $K_2$ in the small interval $[a_i, a_{i+1}]$.
%\begin{prop}\label{prop419}
%If $\sigma\neq 0$ and the L\'evy measure $n(dy)$ satisfies
%(\ref{levy110}), then
%\begin{eqnarray}
%E|K_2(t,a_{i+1})-K_2(t,a_i)|^p\leq c |a_{i+1}-a_i|^{p\over 2}.
%\end{eqnarray}
%\end{prop}
%
%\medskip
%Recall (\ref{levy9}) and denote
%\begin{eqnarray*}
%J^*(X_s,X_{s-},a_i,a_{i+1})&=&J^*(X_s,X_{s}-\Delta X_s,a_i,a_{i+1})\\
%&:=&\Big[J_1(s,a_{i+1})-J_1(s,a_i)+J_2(s,a_{i+1})-J_2(s,a_i)\Big]1_{\{|\Delta
%X_s|<1\}}.
%\end{eqnarray*}
%From the proof of the above proposition, we can easily see that
%\begin{cor}\label{cor4}
%Assume (\ref{levy110}) is satisfied, then for any $p\geq 2$,
%\begin{eqnarray*}
%E\left(\int_0^t\int_{-\infty}^\infty J^*(X_s,X_s- y, a_i, a_{i+1})
%n(dy)ds\right)^p\leq c(t,\sigma,p)|a_{i+1}-a_i|^{{2\over 3}p}.
%\end{eqnarray*}
%\end{cor}
%The estimate will be used in the next section.
\section{The local time rough path}

The $p$-variation ($p>2$) result of the local time enables one to
use Young's integration theory to define $\int_{-\infty}^\infty
g(x)d_x L^x_t $ for $g$ being of bounded $q$-variation when $1\leq
q<2$. This is because in this case, for any $q\in [1,2)$, one can
always find a constant $p>2$ such that the condition ${1\over
p}+{1\over q}>1$ for the existence of the Young integral is
satisfied. However, when $q\geq 2$, Young integral is no longer well
defined. We have to use a new integration theory. Lyons' integration
of rough path provides a way to push the result further. But from
\cite{terry}, generally, we cannot expect to have an integration
theory to define integrals such as $\int_{-\infty}^\infty
g(x)d_xL^x_t $. However, inspired by the method in Chapter 6 in
\cite{terry}, we can treat $Z_x:=(L^x_t,g(x) )$ as a process of
variable $x$ in $R^2$. Then it's easy to know that $Z_x$ is of
bounded $\hat q$-variation in $x$, where $\hat q=q$, if $q>2$, and
$\hat q>2$ can be taken any number when $q=2$. 
%Most of the analysis
%in this section can work for $2\leq q<4$, especially we will
%establish the convergence of smooth rough path in the
%$\theta$-variation topology for any $\theta\in (q, 4)$ so to obtain
%${\bf Z}^1_{a,b}$ and ${\bf Z}^2_{a,b}$. 
In the following, we only consider the case that $2\leq q<3$.  We obtain
the existence of the geometric rough path ${\bf Z}=(1,{\bf Z}^1,{\bf
Z}^2)$ associated to $Z_{.}$.
%If we define a
%one-form  $\hat f(x,y)(v,w):=yv$, then $\int_{-\infty}^\infty
%g(x)d_xL^x_t :=\int_{-\infty}^\infty \hat f(Z)d{\bf Z}^1$.

We assume a sightly stronger condition than (\ref{levy110}) for the L\'evy measure: there
exists a constant $\varepsilon>0$ such that
\begin{eqnarray}\label{levynew36}
\int_{R\setminus \{0\}}(|y|^{1+{1\over
q}-(3-q)\varepsilon}\wedge1)n(dy)<\infty.
\end{eqnarray}
We will prove with this condition, the desired geometric rough path 
${\bf Z}=(1,{\bf Z}^1,{\bf
Z}^2)$ is well defined. We need to point out that in the following when we consider the control function and the convergence of the first level path, condition (\ref{levy110}) is still adequate. But we need (\ref{levynew36}) 
in the convergence of the second level path.
Denote $\delta={1\over q}-{(3-q)}\varepsilon$. Note when
$q=2$, $\delta={1\over 2}-\varepsilon$. So condition (\ref{levynew36}) becomes:
there exists $\varepsilon>0$ such that
\begin{eqnarray}\label{levynew37}
\int_{R\setminus \{0\}}(|y|^{{3\over
2}-\varepsilon}\wedge1)n(dy)<\infty.
\end{eqnarray}
Later in this section, we will see under this condition, the
integral $\int_{-\infty}^\infty L_t^xdL_t^x$ can be well-defined as
a rough path integral. Also note $\inf\limits_{2\leq q<3}
\delta(q,\varepsilon)={1\over 3}$. So under the condition
\begin{eqnarray}\label{levynew38}
\int_{R\setminus \{0\}}(|y|^{{4\over
3}}\wedge1)n(dy)<\infty,
\end{eqnarray}
(\ref{levynew36}) is satisfied for any $2\leq q<3$. 
In this case, our results imply that we can construct the geometric rough path for any
$g$ being of finite $q$-variation, where $2\leq q<3$ can be arbitrary.

 Recall the
$\theta$-variation metric $d_\theta$ on $C_{0,\theta}(\Delta,
T^{([\theta])} (R^2) )$ defined in \cite{terry},
\begin{eqnarray*}
d_\theta({\bf Z},{\bf Y})=\max\limits_{1\leq i\leq
[\theta]}d_{i,\theta}({\bf Z}^i,{\bf Y}^i)=\max\limits_{1\leq i\leq
[\theta]}\sup_D\left(\sum_l|{\bf Z}_{x_{l-1},x_l}^i-{\bf
Y}_{x_{l-1},x_l}^i|^{\theta\over i}\right)^{i\over \theta}.
\end{eqnarray*}

Assume condition (\ref{levy110}) through to Proposition \ref{propo1}.
 Let $[x',x'']$ be any interval in $R$.
From the proof of Theorem \ref{newlem}, for any $p\geq 2$, we know
there exists a constant $c>0$ such that
\begin{eqnarray}\label{control4}
E|L_t^{b}-L_t^{a}|^p\leq c|b-a|^{p\over 2},
\end{eqnarray}
i.e. $L_t^x$ satisfies H$\ddot{\rm o}$lder condition in
\cite{terry} with exponent ${1\over 2}$. First we consider the case
when $g$ is continuous. Recall in \cite{terry}, a control $w$ is a
continuous super-additive function on $\Delta:=\{(a,b):x'\leq
a<b\leq x''\}$ with values in $[0,\infty)$ such that $w(a,a)=0$.
Therefore
\begin{eqnarray*}
w(a,b)+w(b,c)\leq w(a,c),\ \ \ for\ any\ (a,b), (b,c)\in \Delta.
\end{eqnarray*}
If $g(x)$ is of bounded $q$-variation, we can find a control $w$
s.t.
\begin{eqnarray*}
|g(b)-g(a)|^q \leq w(a,b),
\end{eqnarray*}
for any $(a,b)\in \Delta:=\{(a,b):x'\leq a<b\leq x''\}$. It is
obvious that $w_1(a,b):=w(a,b)+(b-a)$ is also a control of $g$. Set
$h={1\over q}$, it is trivial to see for any $\theta>q$ (so
$h\theta>1$) we have,
\begin{eqnarray}\label{may3}
|g(b)-g(a)|^{\theta}\leq w_1(a,b)^{h\theta},\ for \ any \  (a,b)\in
\Delta .
\end{eqnarray}
Considering (\ref{control4}), we can see $Z_x$ satisfies, for such
$h={1\over q}$, and any $\theta>q$, there exists a
constant $c$ such that
\begin{eqnarray}\label{control5}
E|Z_b-Z_a|^{\theta}\leq cw_1(a,b)^{h\theta},\ for \ any \ (a,b)\in
\Delta.
\end{eqnarray}
For any $m\in N$, define a continuous and bounded variation path
$Z(m)$ by
\begin{eqnarray}\label{levy22}
Z(m)_x:=Z_{x_{l-1}^m}+{{w_1(x)-w_1(x_{l-1}^m)}\over
{w_1(x_l^m)-w_1(x_{l-1}^m)}}\Delta_l^m Z,
\end{eqnarray}
if $x_{l-1}^m\leq x<x_l^m$, for $l=1,\cdots,2^m$, and $\Delta_l^m
Z=Z_{x_l^m}-Z_{x_{l-1}^m}$. Here
$D_m:=\{x'=x_0^m<x_1^m<\cdots<x_{2^m}^m=x''\}$ is a partition of
$[x',x'']$ such that
\begin{eqnarray*}
w_1(x_l^m)-w_1(x_{l-1}^m)={1\over {2^m}}w_1(x',x''),
\end{eqnarray*}
where $w_1(x):=w_1(x',x)$. It is obvious that $x_l^m-x_{l-1}^m\leq
{1\over {2^m}}w_1(x',x'')$ and by the superadditivity of the control
function $w_1$,
\begin{eqnarray*}
w_1(x_{l-1}^m,x_l^m)\leq w_1(x_l^m)-w_1(x_{l-1}^m)={1\over
2^m}w_1(x',x'').
\end{eqnarray*}
 The corresponding smooth rough path ${\bf Z}(m)$ is built by
taking its iterated path integrals, i.e. for any $(a,b)\in \Delta$,
\begin{eqnarray}
{\bf
Z}(m)_{a,b}^j=\int_{a<x_1<\cdots<x_j<b}dZ(m)_{x_1}\otimes\cdots\otimes
dZ(m)_{x_j}.
\end{eqnarray}
 In the following, we
will prove $\{{\bf Z}(m)\}_{m\in N}$ converges to a geometric rough
path ${\bf Z}$ in the $\theta$-variation topology when $2\le q<3$.
We call ${\bf Z}$ the canonical geometric rough path associated to
$Z$.
\begin{rmk}
The bounded variation process $Z(m)_x$ is a generalized Wong-Zakai
approximation to the process $Z$ of bounded $\hat q$-variation. Here
we %use the control function $w_1$ to
divide $[x',x'']$ by equally
partitioning the range of $w_1$. We then use (\ref {levy22}) to form
the piecewise curved approximation to $Z$. Note here Wong-Zakai's
standard piecewise linear approximation does not work immediately.
%The standard Wong-Zakai approximation is a special case when
%$w_1(x)=x$, so the approximation $Z(m)$ was piece-wise linear
%(straight line). Here we have curved, but bounded variation
%approximation.
\end{rmk}

Let's first look at the first level path ${\bf Z}(m)_{a,b}^1$. The
method is similar to Chapter 4 in \cite{terry} for Brownian motion.
Similar to Proposition 4.2.1 in \cite{terry}, we can prove that for
all $n\in N$, $m\mapsto\sum\limits_{k=1}^{2^n}|{\bf
Z}(m)_{x_{k-1}^n,x_k^n}^1|^\theta$ is increasing and for $m\geq n$,
\begin{eqnarray}\label{cfeng1}
{\bf Z}(m)^1_{x_{k-1}^n,x_k^n}={\bf
Z}(m+1)^1_{x_{k-1}^n,x_k^n}=Z_{x_k^n}-Z_{x_{k-1}^n}.
\end{eqnarray}
Let ${\bf Z}_{a,b}^1=Z_b-Z_a$. Then (\ref{control5}) implies $E|{\bf
Z}_{a,b}^1|^\theta\leq cw_1(a,b)^{h\theta}$. For such points
$\{x_k^n\}$, $k=1,\cdots,2^n$, $n=1,2,\cdots$, defined above we
still have the inequality in Proposition 4.1.1 in \cite{terry},
\begin{eqnarray}\label{series1}
E\sup_D\sum_{l}|{\bf
Z}_{x_{l-1},x_l}^1|^\theta&\leq&C(\theta,\gamma)E\sum_{n=1}^\infty
n^{\gamma}
\sum_{k=1}^{2^n}|{\bf Z}_{x_{k-1}^n,x_k^n}^1|^\theta\nonumber\\
&\leq&C_1\sum_{n=1}^\infty
n^{\gamma}({1\over{2^n}})^{h\theta-1}w_1(x',x'')^{h\theta},
\end{eqnarray}
for constant $C_1=C(\theta, \gamma)c$. Since $h\theta-1>0$, the
series on the right-hand side of ($\ref{series1}$) is convergent, so
that $\sup_D\sum_l|{\bf Z}_{x_{l-1},x_l}^1|^\theta<\infty$ almost
surely. This shows that ${\bf Z}^1$ has finite $\theta$-variation
almost surely. Moreover, for any $\gamma>\theta-1$, there exists a
constant $C_1(\theta,\gamma,c)>0$ such that
\begin{eqnarray}\label{newww1}
E\sup_m\sup_D\sum_{l}|{\bf Z}(m)_{x_{l-1},x_l}^1|^\theta&\leq&
C(\theta,\gamma)E\sup_m\sum_{n=1}^\infty n^{\gamma}
\sum_{k=1}^{2^n}|{\bf Z}(m)_{x_{k-1}^n,x_k^n}^1|^\theta\nonumber\\
&\leq&C(\theta,\gamma)E\sum_{n=1}^\infty n^{\gamma}
\sum_{k=1}^{2^n}|{\bf Z}_{x_{k-1}^n,x_k^n}^1|^\theta\nonumber\\
&\leq& C_1(\theta,\gamma,c)\sum_{n=1}^\infty
n^{\gamma}({1\over{2^n}})^{h\theta-1}w_1(x',x'')^{h\theta}\\
&<&\infty.\nonumber
\end{eqnarray}
So
\begin{eqnarray}\label{seq1}
\sup_m\sup_D\sum_{l}|{\bf Z}(m)_{x_{l-1},x_l}^1|^\theta< \infty \ \
a.s.
\end{eqnarray}
This means that ${\bf Z}(m)_{a,b}^1$ have finite $\theta$-variation
uniformly in $m$. And furthermore, from (\ref{cfeng1}) and some
standard arguments,
\begin{eqnarray}\label{cfeng2}
E\sum_{n=1}^\infty n^{\gamma} \sum_{k=1}^{2^n}|{\bf
Z}(m)_{x_{k-1}^n,x_k^n}^1-{\bf Z}_{x_{k-1}^n,x_k^n}^1|^\theta\leq
C({1\over{2^m}})^{{h\theta-1}\over 2},
\end{eqnarray}
where $C$ depends on $\theta$, $h$, $w_1(x',x'')$, and $c$ in
(\ref{control5}). By Proposition 4.1.2 in \cite{terry}, Jensen's
inequality and (\ref{cfeng2}),
\begin{eqnarray}\label{control6}
E\sum_{m=1}^\infty\sup_D\Big(\sum_l|{\bf
Z}(m)_{x_{l-1},x_l}^1-{\bf
Z}_{x_{l-1},x_l}^1|^\theta\Big)^{1\over\theta}
&\leq&E\sum_{m=1}^\infty\left(\sum_{n=1}^\infty n^{\gamma}
\sum_{k=1}^{2^n}|{\bf Z}(m)_{x_{k-1}^n,x_k^n}^1-{\bf Z}_{x_{k-1}^n,x_k^n}^1|^\theta\right)^{1\over \theta}\nonumber\\
&\leq&\sum_{m=1}^\infty\left(E\sum_{n=1}^\infty n^{\gamma}
\sum_{k=1}^{2^n}|{\bf Z}(m)_{x_{k-1}^n,x_k^n}^1-{\bf Z}_{x_{k-1}^n,x_k^n}^1|^\theta\right)^{1\over \theta}\nonumber\\
&\leq&C\sum_{m=1}^\infty({1\over{2^m}})^{{h\theta-1}\over
{2\theta}}\nonumber\\
&<&\infty,
\end{eqnarray}
for $h\theta>1$. So we obtain
\begin{thm}\label{theor1}
Let $L_t^x$ be the local time of the time homogeneous L\'evy process
$X_t$ given by (\ref{tanakaf1}), and $g$ be a continuous function of
bounded $q$-variation. Assume $q\geq 1$, $\sigma\neq 0$ and the
L\'evy measure $n(dy)$ satisfies (\ref{levy110}). Then for any
$\theta>q$, the continuous process $Z_x=(L_t^x,g(x))$ satisfying
(\ref{control5}), we have
\begin{eqnarray}
\sum_{m=1}^\infty\sup_D\Big(\sum_l|{\bf Z}(m)_{x_{l-1},x_l}^1-{\bf
Z}_{x_{l-1},x_l}^1|^\theta\Big)^{1\over\theta}<\infty\ a.s..
\end{eqnarray}
In particular, $({\bf Z}(m)_{a,b}^1)$ converges to $({\bf
Z}_{a,b}^1)$ in the $\theta$-variation distance a.s. for any
$(a,b)\in \Delta$.
\end{thm}

We next consider the second level path ${\bf Z}(m)_{a,b}^2$. As in
\cite{terry}, we can also see that if $m\leq n$,
\begin{eqnarray}\label{cfeng3}
{\bf Z}(m)_{x_{k-1}^n,x_k^n}^2=2^{2(m-n)-1}(\Delta_l^mZ)^{\otimes2},
\end{eqnarray}
where $l$ is chosen such that $x_{l-1}^m\le x_{k-1}^n<x_k^n\leq
x_l^m$; if $m>n$,
$${\bf Z}(m)_{x_{k-1}^n,x_k^n}^2={1\over 2}\Delta_k^n
Z\otimes\Delta_k^n Z+{1\over 2}\sum_{\stackrel
{r,l=2^{m-n}(k-1)+1}{r<l}}^{2^{m-n}k}(\Delta_r^mZ\otimes\Delta_l^mZ-\Delta_l^mZ\otimes\Delta_r^mZ),$$
so
\begin{eqnarray}\label{cfeng11}
&&{\bf Z}(m+1)_{x_{k-1}^n,x_k^n}^2-{\bf
Z}(m)_{x_{k-1}^n,x_k^n}^2\nonumber\\
&=&{1\over
2}\sum_{l=2^{m-n}(k-1)+1}^{2^{m-n}k}(\Delta_{2l-1}^{m+1}Z\otimes\Delta_{2l}^{m+1}Z-\Delta_{2l}^{m+1}Z\otimes\Delta_{2l-1}^{m+1}Z),
\end{eqnarray}
$k=1,\cdots,2^n$. Similar to the proof of Proposition 4.3.3 in
\cite{terry}, we have
\begin{prop}\label{propo1}
Assume $g$ is a continuous function of finite $q$-variation with a real number $q\geq 2$,  and the 
L\'evy measure satisfies (\ref{levy110}). Let $\theta>q$. Then for $m\leq n$,
\begin{eqnarray}\label{may12}
\sum_{k=1}^{2^n}E|{\bf Z}(m+1)_{x_{k-1}^n,x_k^n}^2-{\bf
Z}(m)_{x_{k-1}^n,x_k^n}^2|^{\theta\over 2}\leq
C({1\over{2^{n+m}}})^{{\theta h-1}\over 2},
\end{eqnarray}
where $C$ depends on $\theta$, $h(:={1\over q})$, $w_1(x',x'')$, and
$c$ in (\ref{control5}).
\end{prop}

The main step to establish the geometric rough path integral over
$Z$ is the following estimate. The proof is rather complicated. We will use Lemma \ref{lem8} 
about the correlation of $K_2(t,a_{i+1})-K_2(t,a_{i})$ and
$K_2(t,a_{j+1})-K_2(t,a_{j})$, and (\ref{2june2009}) for the term $K_3$, 
for $\xi={{q-2}\over {2q}}+(3-q)\varepsilon$.

\begin{prop}\label{propo2}
Assume $g$ is a continuous function of finite $q$-variation with a real number $q\in [2,3)$,  and the 
L\'evy measure satisfies (\ref{levynew36}). Let $q<\theta<3$. Then for $m>n$, we have
\begin{eqnarray}\label{may13}
E|{\bf Z}(m+1)_{x_{k-1}^n,x_k^n}^2-{\bf
Z}(m)_{x_{k-1}^n,x_k^n}^2|^{\theta\over 2} \leq C
\Big[({1\over{2^n}})^{\theta\over 4}({1\over{2^m}})^{{1\over
2}h\theta}+({1\over{2^n}})^{\theta\over
2}({1\over{2^m}})^{{{3-q}\over 2}{\varepsilon\theta}}\Big],
\end{eqnarray}
where $C$ is a generic constant and also depends on $\theta$,
$h(:={1\over q})$, $w_1(x',x'')$, and $c$ in (\ref{control5}).
\end{prop}

\begin{lem}\label{lem8}
Assume the L\'evy measure satisfies (\ref{levynew34}) with
$0\leq \xi\leq {1\over 6}$, then for any $a_0<a_1<\cdots<a_m$,
%\begin{eqnarray}\label{levy41}
%&&\Big|E\int_0^{t+}\int_{|y|<1}1_{\{a_i< X_{s-}\leq a_{i+1}\}}y
% \tilde N_p(dsdy)\big(K_2(t,a_{j+1})-K_2(t,a_j)\big)\Big|\nonumber\\
% &\leq& \cases {c(t,\sigma) (a_{i+1}-a_i), \ \quad \quad \quad \quad \quad
% \quad \quad \quad {\rm \ when }
%\ 0\leq i=j\leq m,\cr c(t,\sigma)[(a_{i+1}-a_i)^{4\over
%3}+(a_{j+1}-a_j)^{4\over 3}],\ \ \ {\rm \ when } \ 0\leq i\neq j\leq
%m,}
%\end{eqnarray}
%and
\begin{eqnarray}\label{levy42}
&&\Big|E\big(K_2(t,a_{i+1})-K_2(t,a_i)\big)\big(K_2(t,a_{j+1})-K_2(t,a_j)\big)\Big|\nonumber\\
 &\leq& \cases {c(t,\sigma) (a_{i+1}-a_i), \ \ \ \ \ \ \ \ \ \ \ \ \quad \quad \quad \quad  \quad \quad \quad {\rm \ when }
\ 0\leq i=j\leq m,\cr
c(t,\sigma)[(a_{i+1}-a_i)^{1+2\xi}+(a_{j+1}-a_j)^{1+2\xi}],\ \ \
{\rm \ when } \ 0\leq i\neq j\leq m.}
\end{eqnarray}
\end{lem}
{\bf Proof}:
When $i=j$, (\ref{levy42}) follows from (\ref{levy10})
directly. %it is easy to obtain that
%\begin{eqnarray}
%E\big(K_2(t,a_{i+1})-K_2(t,a_{i})\big)^2 \leq  c(t) |a_{i+1}-a_{i}|.
%\end{eqnarray}
%When $i\neq j$,
Now we consider the case when $i\neq j$. Without losing generality,
we assume that $i<j$. From (\ref{levy9}), it is easy to see that
\begin{eqnarray*}
&&E\big(K_2(t,a_{i+1})-K_2(t,a_{i})\big)\cdot\big(K_2(t,a_{j+1})-K_2(t,a_{j})\big)\\
&=&E\left [\int_0^{t+}\int_R J_*(X_s, X_s-y,a_i,a_{i+1})\tilde
N_p(dyds)\cdot\int_0^{t+}\int_R J_*(X_s, X_s-y,a_j,a_{j+1})\tilde
N_p(dyds)\right ]\\
 &=&E\int_0^{t}\int_R J^*(X_s, X_s-y,a_i,a_{i+1})J^*(X_s,
X_s-y,a_j,a_{j+1})n(dy)ds.
%&&+E\int_0^{t}\int_R J^*(X_s,
%X_s-y,a_i,a_{i+1})n(dy)ds\cdot\int_0^{t}\int_R J^*(X_s,
%X_s-y,a_j,a_{j+1})n(dy)ds.
\end{eqnarray*}
%By using the Cauchy-Schwarz inequality and Corollary \ref{cor4}, we
%know,
%\begin{eqnarray*}
%&&\left|E\int_0^{t}\int_R J^*(X_s,
%X_s-y,a_i,a_{i+1})n(dy)ds\cdot\int_0^{t}\int_R J^*(X_s,
%X_s-y,a_j,a_{j+1})n(dy)ds\right|\\
%&\leq& \left(E\Big(\int_0^{t}\int_R J^*(X_s, X_s-y,a_i,a_{i+1})
%n(dy)ds\Big)^2\right)^{1\over 2} \left(E\Big(\int_0^{t}\int_R
%J^*(X_s,
%X_s-y,a_j,a_{j+1})n(dy)ds\Big)^2\right)^{1\over 2}\\
%&\leq & c(t,\sigma)(a_{i+1}-a_i)^{2\over 3}(a_{j+1}-a_j)^{2\over 3}.
%\end{eqnarray*}
But from (\ref{levy57}), we know
\begin{eqnarray*}
\int_0^{t}\int_R J_*(X_s, X_s-y,a_i,a_{i+1})J_*(X_s,
X_s-y,a_j,a_{j+1})n(dy)ds=A_1+A_2+A_3+A_4,
\end{eqnarray*}
where
\begin{eqnarray*}
A_1&:=&\int_0^{t}\int_R\Big[-(a_{i+1}-a_i)(a_j-X_{s-})\Big]1_{\{X_s\leq
a_i\}}1_{\{a_j<X_{s-}\leq
a_{j+1}\}}1_{\{|y|<1\}}n(dy)ds\\
&&+\int_0^{t}\int_R\Big[(a_{i+1}-a_i)(a_{j+1}-a_j)\Big]1_{\{X_s\leq
a_i\}}1_{\{X_{s-}> a_{j+1}\}}1_{\{|y|<1\}}n(dy)ds,\\
A_2&:=&\int_0^{t}\int_R\Big[-(a_{i+1}-X_s)(a_j-X_{s-})\Big]1_{\{a_i<
X_s\leq
a_{i+1}\}}1_{\{a_j<X_{s-}\leq a_{j+1}\}}1_{\{|y|<1\}}n(dy)ds\\
&&+ \int_0^{t}\int_R (a_{i+1}-X_s)(a_{j+1}-a_j)1_{\{a_i< X_s\leq
a_{i+1}\}}1_{\{X_{s-}> a_{j+1}\}}1_{\{|y|<1\}}n(dy)ds,\\
A_3&:=&\int_0^{t}\int_R\Big[(a_{i+1}-a_i)(X_s-a_j)\Big]1_{\{X_{s-}\leq
a_i\}}1_{\{a_j<X_{s}\leq
a_{j+1}\}}1_{\{|y|<1\}}n(dy)ds\\
&&+\int_0^{t}\int_R(a_{i+1}-a_i)(a_{j+1}-a_j)1_{\{X_{s-}\leq
a_i\}}1_{\{X_{s}> a_{j+1}\}}1_{\{|y|<1\}}n(dy)ds,\\
A_4&:=&\int_0^{t}\int_R\Big[-(X_{s-}-a_{i+1})(X_s-a_j)1_{\{a_i<
X_{s-}\leq a_{i+1}\}}1_{\{a_j<X_{s}\leq
a_{j+1}\}}1_{\{|y|<1\}}n(dy)ds\\
 &&+\int_0^{t}\int_R\Big[-(X_{s-}-a_{i+1})(a_{j+1}-a_j)\Big]1_{\{a_i<
X_{s-}\leq a_{i+1}\}}1_{\{X_{s}>
a_{j+1}\}}\Big]1_{\{|y|<1\}}n(dy)ds.
\end{eqnarray*}
To estimate $|EA_1|$, we notice that
\begin{eqnarray*}
|E A_1 |&\leq &E\int_0^t\int_{-\infty}^\infty
(a_{i+1}-a_i)(a_j-X_{s-})1_{\{X_s\leq
a_i\}}1_{\{a_j<X_{s}-y\leq a_{j+1}\}}1_{\{|y|< 1\}} n(dy)ds\\
%&&+E\Bigg[\int_0^t\int_{-\infty}^\infty (a_{i+1}-a_i)1_{\{X_s\leq
%a_i\}}1_{\{a_j<X_{s}-y\leq a_{j+1}\}}1_{\{|y|< 1\}} n(dy)ds\\
%&&\hskip 1cm\cdot\int_0^t\int_{-\infty}^\infty (a_j-X_s)1_{\{X_s\leq
%a_i\}}1_{\{a_j<X_{s}-y\leq a_{j+1}\}}1_{\{|y|< 1\}} n(dy)ds\Bigg]\\
%&&+E\Bigg[\int_0^t\int_{-\infty}^\infty (a_{i+1}-a_i)1_{\{X_s\leq
%a_i\}}1_{\{a_j<X_{s}-y\leq a_{j+1}\}}1_{\{|y|< 1\}} n(dy)ds\\
%&&\hskip 1cm\cdot\int_0^t\int_{-\infty}^\infty
%(a_{j+1}-a_j)1_{\{X_s\leq
%a_i\}}1_{\{X_{s}-y> a_{j+1}\}}1_{\{|y|< 1\}} n(dy)ds\Bigg]\\
&&+E\int_0^t\int_{-\infty}^\infty(a_{i+1}-a_i)
(a_{j+1}-a_j)1_{\{X_s\leq a_i\}}1_{\{X_{s}-y> a_{j+1}\}}1_{\{|y|<
1\}} n(dy)ds\\
%&&+E\Bigg[\int_0^t\int_{-\infty}^\infty (a_{i+1}-a_i)1_{\{X_s\leq
%a_i\}}1_{\{X_{s}-y> a_{j+1}\}}1_{\{|y|<
%1\}} n(dy)ds\\
%&&\hskip 1cm\cdot\int_0^t\int_{-\infty}^\infty(a_{j+1}-a_j)
%1_{\{X_s\leq a_i\}}1_{\{X_{s}-y> a_{j+1}\}}1_{\{|y|<
%1\}} n(dy)ds\Bigg]\\
%&&+E\Bigg[\int_0^t\int_{-\infty}^\infty (a_{i+1}-a_i)1_{\{X_s\leq
%a_i\}}1_{\{X_{s}-y> a_{j+1}\}}1_{\{|y|<
%1\}} n(dy)ds\\
%&&\hskip 1cm\cdot\int_0^t\int_{-\infty}^\infty (a_j-X_s)1_{\{X_s\leq
%a_i\}}1_{\{a_j<X_{s}-y\leq a_{j+1}\}}1_{\{|y|< 1\}} n(dy)ds\Bigg]\\
&:=&A_{11}+A_{12}.
\end{eqnarray*}
Let's estimate every term on the righthand side of the above
inequality. By the occupation times formula, Fubini theorem,
Jensen's inequality, similar as before, we have :
\begin{eqnarray*}
A_{11}
%&:=&E\int_0^t\int_{-\infty}^\infty
%(a_{i+1}-a_i)(a_j-X_s)1_{\{X_s\leq
%a_i\}}1_{\{a_j<X_{s}-y\leq a_{j+1}\}}1_{\{|y|< 1\}} n(dy)ds\\
&= &{1\over {\sigma^2}}E
\int_{-\infty}^{a_i}L_t^x\int_{x-a_{j+1}}^{x-a_j}
(a_{i+1}-a_i)(x-a_j-y)1_{\{|y|< 1\}} n(dy)dx\\
&\leq &{1\over {\sigma^2}}(a_{i+1}-a_i)^{{1\over 2}+\xi}E\Bigg[
\int_{-\infty}^{a_i-a_{j+1}}|y|^{{1\over
2}-\xi}\int^{y+a_{j+1}}_{y+a_j}L_t^x
|y|1_{\{|y|< 1\}}dx n(dy)\\
&&\hskip 2cm+\int^{a_i-a_j}_{a_i-a_{j+1}}|y|^{{1\over
2}-\xi}\int^{a_i}_{y+a_j}L_t^x
|y|1_{\{|y|< 1\}}dx n(dy)\Bigg]\\
&\leq &{1\over {\sigma^2}}(a_{i+1}-a_i)^{{1\over 2}+\xi}(a_{j+1}-a_j)(\sup_x
EL_t^x)\int_{-\infty}^{a_i-a_j}|y|^{{3\over 2}-\xi}1_{\{|y|< 1\}}n(dy)\\
&\leq & c(t,\sigma)(a_{i+1}-a_i)^{{1\over 2}+\xi}(a_{j+1}-a_j).
\end{eqnarray*}
In the same way, we can have
\begin{eqnarray*}
A_{12}&\leq &{1\over {\sigma^2}}(a_{i+1}-a_i)^{{1\over
2}+\xi}(a_{j+1}-a_j)\sup_x
E(L_t^x)\int_{-\infty}^{a_i-a_{j+1}}|y|^{{3\over 2}-\xi}|y|1_{\{|y|<
1\}}n(dy)\\
&\leq& c(t,\sigma)(a_{i+1}-a_i)^{{1\over 2}+\xi}(a_{j+1}-a_j).
%&&A_{15}\leq {2\over {\sigma^2}}(a_{i+1}-a_i)(a_{j+1}-a_j)\sup_x
%E(L_t^x)^2\int_{-\infty}^{a_i-a_{j+1}}|y|1_{\{|y|< 1\}}n(dy)\leq
%c(t,\sigma)(a_{i+1}-a_i)(a_{j+1}-a_j)\\
%&&A_{16}\leq {2\over {\sigma^2}}(a_{i+1}-a_i)(a_{j+1}-a_j)\sup_x
%E(L_t^x)\int_{-\infty}^{a_i-a_{j}}|y|1_{\{|y|< 1\}}n(dy)\leq
%c(t,\sigma)(a_{i+1}-a_i)(a_{j+1}-a_j).
\end{eqnarray*}
Therefore, we get
\begin{eqnarray*}
|EA_1| \leq
c(t,\sigma)((a_{i+1}-a_i)^{1+2\xi}+(a_{j+1}-a_j)^{1+2\xi}).
\end{eqnarray*}
Using the same method, we can have a similar estimation in the
other cases:
\begin{eqnarray*}
|E A_2 | &\leq &E\int_0^t\int_{-\infty}^\infty
(a_{i+1}-X_s)(X_s-y-a_j)1_{\{a_i<X_s\leq
a_{i+1}\}}1_{\{a_j<X_{s}-y\leq a_{j+1}\}}1_{\{|y|< 1\}} n(dy)ds\\
%&&+E\Bigg[\int_0^t\int_{-\infty}^\infty
%(a_{i+1}-X_s)1_{\{a_i<X_s\leq
%a_{i+1}\}}1_{\{a_j<X_{s}-y\leq a_{j+1}\}}1_{\{|y|< 1\}} n(dy)ds\\
%&&\hskip 1cm \cdot \int_0^t\int_{-\infty}^\infty
%(a_j-X_s)1_{\{a_i<X_s\leq
%a_{i+1}\}}1_{\{a_j<X_{s}-y\leq a_{j+1}\}}1_{\{|y|< 1\}} n(dy)ds\Bigg]\\
%&&+E\Bigg[\int_0^t\int_{-\infty}^\infty
%(a_{i+1}-X_s)1_{\{a_i<X_s\leq
%a_{i+1}\}}1_{\{a_j<X_{s}-y\leq a_{j+1}\}}1_{\{|y|< 1\}} n(dy)ds\\
%&&\hskip 1cm\cdot \int_0^t\int_{-\infty}^\infty
%(a_{j+1}-a_j)1_{\{a_i<X_s\leq
%a_{i+1}\}}1_{\{X_{s}-y> a_{j+1}\}}1_{\{|y|< 1\}} n(dy)ds\Bigg]\\
&&+E\int_0^t\int_{-\infty}^\infty
(a_{i+1}-X_s)(a_{j+1}-a_j)1_{\{a_i<X_s\leq a_{i+1}\}}1_{\{X_{s}-y>
a_{j+1}\}}1_{\{|y|< 1\}} n(dy)ds\\
%&&+E\Bigg[\int_0^t\int_{-\infty}^\infty
%(a_{i+1}-X_s)1_{\{a_i<X_s\leq a_{i+1}\}}1_{\{X_{s}-y>
%a_{j+1}\}}1_{\{|y|< 1\}} n(dy)ds\\
%&&\hskip 1cm\cdot\int_0^t\int_{-\infty}^\infty
%(a_{j+1}-a_j)1_{\{a_i<X_s\leq a_{i+1}\}}1_{\{X_{s}-y>
%a_{j+1}\}}1_{\{|y|< 1\}} n(dy)ds\Bigg]\\
%&&+E\Bigg[\int_0^t\int_{-\infty}^\infty
%(a_{i+1}-X_s)1_{\{a_i<X_s\leq a_{i+1}\}}1_{\{X_{s}-y>
%a_{j+1}\}}1_{\{|y|< 1\}} n(dy)ds\\
%&&\hskip 1cm\cdot\int_0^t\int_{-\infty}^\infty
%(a_j-X_s)1_{\{a_i<X_s\leq a_{i+1}\}}1_{\{a_j<X_{s}-y\leq
%a_{j+1}\}}1_{\{|y|< 1\}}
%n(dy)ds\Bigg]\\
&\leq& c(t,\sigma)((a_{i+1}-a_i)^{1+2\xi}+(a_{j+1}-a_j)^{1+2\xi});
\end{eqnarray*}
and
\begin{eqnarray*}
\hskip -1cm |E A_3| &\leq&E\int_0^t\int_{-\infty}^\infty
(a_{i+1}-a_i)(X_s-a_j)1_{\{X_s-y\leq
a_i\}}1_{\{a_j<X_{s}\leq a_{j+1}\}}1_{\{|y|< 1\}} n(dy)ds\\
%&&+E\Bigg[\int_0^t\int_{-\infty}^\infty (a_{i+1}-a_i)1_{\{X_s-y\leq
%a_i\}}1_{\{a_j<X_{s}\leq a_{j+1}\}}1_{\{|y|< 1\}} n(dy)ds\\
%&&\hskip 1cm\cdot\int_0^t\int_{-\infty}^\infty
%(X_s-a_j)1_{\{X_s-y\leq
%a_i\}}1_{\{a_j<X_{s}\leq a_{j+1}\}}1_{\{|y|< 1\}} n(dy)ds\Bigg]\\
%&&+E\Bigg[\int_0^t\int_{-\infty}^\infty (a_{i+1}-a_i)1_{\{X_s-y\leq
%a_i\}}1_{\{a_j<X_{s}\leq a_{j+1}\}}1_{\{|y|< 1\}} n(dy)ds\\
%&&\hskip 1cm\cdot\int_0^t\int_{-\infty}^\infty
%(a_{j+1}-a_j)1_{\{X_s-y\leq
%a_i\}}1_{\{X_{s}> a_{j+1}\}}1_{\{|y|< 1\}} n(dy)ds\Bigg]\\
&&+E\int_0^t\int_{-\infty}^\infty(a_{i+1}-a_i)
(a_{j+1}-a_j)1_{\{X_s-y\leq a_i\}}1_{\{X_{s}> a_{j+1}\}}1_{\{|y|<
1\}} n(dy)ds\\
%&&+E\Bigg[\int_0^t\int_{-\infty}^\infty (a_{i+1}-a_i)1_{\{X_s-y\leq
%a_i\}}1_{\{X_{s}> a_{j+1}\}}1_{\{|y|<
%1\}} n(dy)ds\\
%&&\hskip 1cm\cdot\int_0^t\int_{-\infty}^\infty(a_{j+1}-a_j)
%1_{\{X_s-y\leq a_i\}}1_{\{X_{s}> a_{j+1}\}}1_{\{|y|<
%1\}} n(dy)ds\Bigg]\\
%&&+E\Bigg[\int_0^t\int_{-\infty}^\infty (a_{i+1}-a_i)1_{\{X_s-y\leq
%a_i\}}1_{\{X_{s}> a_{j+1}\}}1_{\{|y|<
%1\}} n(dy)ds\\
%&&\hskip 1cm\cdot\int_0^t\int_{-\infty}^\infty
%(X_s-a_j)1_{\{X_s-y\leq
%a_i\}}1_{\{a_j<X_{s}\leq a_{j+1}\}}1_{\{|y|< 1\}} n(dy)ds\Bigg]\\
&\leq& c(t,\sigma)((a_{i+1}-a_i)^{1+2\xi}+(a_{j+1}-a_j)^{1+2\xi});
\end{eqnarray*}
and
\begin{eqnarray*}
 |E A_4| &\leq
&E\int_0^t\int_{-\infty}^\infty
(a_{i+1}-X_s+y)(X_s-a_j)1_{\{a_i<X_s-y\leq
a_{i+1}\}}1_{\{a_j<X_{s}\leq a_{j+1}\}}1_{\{|y|< 1\}} n(dy)ds\\
%&&+E\Bigg[\int_0^t\int_{-\infty}^\infty
%(X_s-a_{i+1})1_{\{a_i<X_s-y\leq
%a_{i+1}\}}1_{\{a_j<X_{s}\leq a_{j+1}\}}1_{\{|y|< 1\}} n(dy)ds\\
%&&\hskip 1cm\cdot \int_0^t\int_{-\infty}^\infty
%(X_s-a_j)1_{\{a_i<X_s-y\leq
%a_{i+1}\}}1_{\{a_j<X_{s}\leq a_{j+1}\}}1_{\{|y|< 1\}} n(dy)ds\Bigg]\\
%&&+E\Bigg[\int_0^t\int_{-\infty}^\infty
%(X_s-a_{i+1})1_{\{a_i<X_s-y\leq
%a_{i+1}\}}1_{\{a_j<X_{s}\leq a_{j+1}\}}1_{\{|y|< 1\}} n(dy)ds\\
%&&\hskip 1cm\cdot \int_0^t\int_{-\infty}^\infty
%(a_{j+1}-a_j)1_{\{a_i<X_s-y\leq
%a_{i+1}\}}1_{\{X_{s}> a_{j+1}\}}1_{\{|y|< 1\}} n(dy)ds\Bigg]\\
&&+E\int_0^t\int_{-\infty}^\infty
(a_{i+1}-X_s+y)(a_{j+1}-a_j)1_{\{a_i<X_s-y\leq a_{i+1}\}}1_{\{X_{s}>
a_{j+1}\}}1_{\{|y|< 1\}} n(dy)ds\\
%&&+E\Bigg[\int_0^t\int_{-\infty}^\infty
%(X_s-a_{i+1})1_{\{a_i<X_s-y\leq a_{i+1}\}}1_{\{X_{s}>
%a_{j+1}\}}1_{\{|y|< 1\}} n(dy)ds\\
%&&\hskip 1cm\cdot\int_0^t\int_{-\infty}^\infty
%(a_{j+1}-a_j)1_{\{a_i<X_s-y\leq a_{i+1}\}}1_{\{X_{s}>
%a_{j+1}\}}1_{\{|y|< 1\}} n(dy)ds\Bigg]\\
%&&+E\Bigg[\int_0^t\int_{-\infty}^\infty
%(X_s-a_{i+1})1_{\{a_i<X_s-y\leq a_{i+1}\}}1_{\{X_{s}>
%a_{j+1}\}}1_{\{|y|< 1\}} n(dy)ds\\
%&&\hskip 1cm\cdot\int_0^t\int_{-\infty}^\infty
%(X_s-a_j)1_{\{a_i<X_s-y\leq a_{i+1}\}}1_{\{a_j<X_{s}\leq
%a_{j+1}\}}1_{\{|y|< 1\}}
%n(dy)ds\Bigg]\\
&\leq& c(t,\sigma)((a_{i+1}-a_i)^{1+2\xi}+(a_{j+1}-a_j)^{1+2\xi}).
\end{eqnarray*}
So when $i\neq j$,
\begin{eqnarray*}
|E\big(K_2(t,a_{i+1})-K_2(t,a_{i})\big)\big(K_2(t,a_{j+1})-K_2(t,a_{j})\big)|
\leq c(t,\sigma)((a_{i+1}-a_i)^{1+2\xi}+(a_{j+1}-a_j)^{1+2\xi}).
\end{eqnarray*}
Therefore we proved (\ref{levy42}).
 $\hfill\diamond$
 \vskip5pt
 {\bf
Proof of Proposition \ref{propo2}:} For $m>n$, it is easy to see
that
\begin{eqnarray}\label{control7}
&&E|{\bf Z}(m+1)_{x_{k-1}^n,x_k^n}^2-{\bf Z}(m)_{x_{k-1}^n,x_k^n}^2|^2\nonumber\\
&=&{1\over
4}E\Big|\sum_{l=2^{m-n}(k-1)+1}^{2^{m-n}k}(\Delta_{2l-1}^{m+1}Z\otimes\Delta_{2l}^{m+1}Z-\Delta_{2l}^{m+1}Z\otimes\Delta_{2l-1}^{m+1}Z)\Big|^2\nonumber\\
&=&{1\over 4}\sum_{\stackrel {i,j= 1}{i\neq
j}}^2E\sum_{l,r=2^{m-n}(k-1)+1}^{2^{m-n}k}
(\Delta_{2l-1}^{m+1}Z^i\Delta_{2l}^{m+1}Z^j-\Delta_{2l}^{m+1}Z^i\Delta_{2l-1}^{m+1}Z^j)\nonumber\\
&&\hskip 4cm\cdot(\Delta_{2r-1}^{m+1}Z^i\Delta_{2r}^{m+1}Z^j-\Delta_{2r}^{m+1}Z^i\Delta_{2r-1}^{m+1}Z^j)\nonumber\\
%&=&{1\over 4}\sum_{\stackrel {i,j= 1}{i\neq
%j}}^2\sum_{l,r=2^{m-n}(k-1)+1}^{2^{(m-
%n)}k}\Big[E(\Delta_{2l-1}^{m+1}Z^i\Delta_{2r-1}^{m+1}Z^i)
%E(\Delta_{2l}^{m+1}Z^j\Delta_{2r}^{m+1}Z^j)\nonumber\\
%&&\hskip 4cm
%+E(\Delta_{2l}^{m+1}Z^i\Delta_{2r}^{m+1}Z^i)E(\Delta_{2l-1}^{m+1}Z^j\Delta_{2r
%-1}^{m+1}Z^j)\Big]\nonumber\\
%&&-{1\over 4}\sum_{\stackrel {i,j= 1}{i\neq
%j}}^2\sum_{l,r=2^{m-n}(k-1)+1}^{2^{(m-
%n)}k}\Big[E(\Delta_{2l}^{m+1}Z^i\Delta_{2r-1}^{m+1}Z^i) E(\Delta_{2l
%-1}^{m+1}Z^j\Delta_{2r}^{m+1}Z^j)\nonumber\\
%&&\hskip 4cm
%+E(\Delta_{2l-1}^{m+1}Z^i\Delta_{2r}^{m+1}Z^i)E(\Delta_{2l}^{m+1}Z^j\Delta_{2r
%-1}^{m+1}Z^j)\Big]\nonumber\\
&=&{1\over 4}\sum_{l,r}\Big[E(\Delta_{2l-1}^{m+1}
L_t^x\Delta_{2r-1}^{m+1}
L_t^x)(\Delta_{2l}^{m+1}g(x)\Delta_{2r}^{m+1}g(x))\nonumber\\
&&\hskip1cm+(\Delta_{2l-1}^{m+1}g(x)\Delta_{2r-1}^{m+1}g(x))E(\Delta_{2l}^{m+1}
L_t^x\Delta_{2r}^{m+1}
L_t^x)\Big]\nonumber\\
&&-{1\over 4}\sum_{l,r}\Big[E(\Delta_{2l-1}^{m+1}
L_t^x\Delta_{2r}^{m+1}
L_t^x)(\Delta_{2l}^{m+1}g(x)\Delta_{2r-1}^{m+1}g(x))\nonumber\\
&&\hskip
1cm+(\Delta_{2l-1}^{m+1}g(x)\Delta_{2r}^{m+1}g(x))E(\Delta_{2l}^{m+1}
L_t^x\Delta_{2r-1}^{m+1}L_t^x)\Big]\nonumber\\
 &&-{1\over 4}\sum_{l,r}\Big[E(\Delta_{2l}^{m+1}
L_t^x\Delta_{2r-1}^{m+1}
L_t^x)(\Delta_{2l-1}^{m+1}g(x)\Delta_{2r}^{m+1}g(x))\nonumber\\
&&\hskip 1cm+(\Delta_{2l}^{m+1}g(x)\Delta_{2r-1}^{m+1}g(x))E(\Delta_{2l-1}^{m+1}L_t^x\Delta_{2r}^{m+1}L_t^x)\Big]\nonumber\\
&&+{1\over 4}\sum_{l,r}\Big[E(\Delta_{2l}^{m+1}
L_t^x\Delta_{2r}^{m+1}
L_t^x)(\Delta_{2l-1}^{m+1}g(x)\Delta_{2r-1}^{m+1}g(x))\nonumber\\
&&\hskip
1cm+(\Delta_{2l}^{m+1}g(x)\Delta_{2r}^{m+1}g(x))E(\Delta_{2l-1}^{m+1}
L_t^x\Delta_{2r-1}^{m+1}L_t^x)\Big].
\end{eqnarray}
%From Tanaka's formula (c.f. \cite{ks} and \cite{yor}) and
%\begin{eqnarray*}
%\sum \limits _{x^*_k\le x}\hat L_t^{x^*_k}=\int _0^t1_{\{X(s)\le
%x\}}dV_s=V_t-V_0-\int _0^t 1_{\{X(s)>x\}}dV_s,
%\end{eqnarray*}
%It's easy to know that
%\begin{eqnarray}
%L_t(a)&=&(X_t-a)^+-(X_0-a)^+-\int_0^t 1_{\{X_{s-}>
%a\}}d(bB_s+\tilde
%M_s)\nonumber\\
%&&+\sum_{0\leq s\leq t}[(X_{s-}-a)^+-(X_s-a)^+]1_{\{|\Delta
%X_s|\geq 1\}}\nonumber\\
%&&+\sum_{0\leq s\leq t}[(X_{s-}-a)^+-(X_s-a)^++1_{\{X_{s-}>
%a\}}\Delta X_s]1_{\{|\Delta X_s|< 1\}}\nonumber\\
%&:=&\varphi_t(a)-b\hat B^a_t-\hat M^a_t+K_1(t,a)+K_2(t,a),
%\end{eqnarray}
%\begin{eqnarray*}
%\tilde
%L_t^x&=&(X_t-x)^+-(X_0-x)^+-\int_0^t1_{\{X_s>x\}}dM_s-(V_t-V_0)\\
%&:=&\phi(x)-\int_0^t1_{\{X_s>x\}}dM_s-(V_t-V_0).
%\end{eqnarray*}
%And using the following estimates in the proof of Lemma 2.1 in
%\cite{Zhao33}: for any $\gamma\geq 1$ and $a_i<a_{i+1}$
%\begin{eqnarray*}
%&&| \varphi_t(a_{i+1})- \varphi_t(a_i)|^\gamma\leq 2^\gamma
%(a_{i+1}-a_i)^\gamma,\ E|b\int_0^t1_{\{a_{i}< X_s\leq
%a_{i+1}\}}ds|^\gamma\leq c(t,b,\gamma)(a_{i+1}-a_i)^{\gamma}\\
%&&E|\sigma\int_0^t1_{\{a_{i}< X_s\leq a_{i+1}\}}dB_s|^\gamma\leq
%c(t,\sigma,\gamma)(a_{i+1}-a_i)^{\gamma\over 2}.
%\end{eqnarray*}
The main difficulty is to estimate the following expectation which
can be derived from Tanaka's formula:
\begin{eqnarray*}
 &&E\Big[\Delta_{2r-1}^{m+1}
L_t^x\Delta_{2l-1}^{m+1}L_t^x\Big]\\
&=&E\Big[\big(L_t(x_{2r-1}^{m+1})-
L_t(x_{2r-2}^{m+1})\big)\big(L_t(x_{2l-1}^{m+1})-
L_t(x_{2l-2}^{m+1})\big)\Big]\\
&=&E\Big[\varphi_t(x_{2r-1}^{m+1})-\varphi_t(x_{2r-2}^{m+1})-b\int_0^t1_{\{x_{2r-2}^{m+1}<
X_{s-}\leq x_{2r-1}^{m+1}\}}ds \\
&&\hskip 0.5cm-\sigma\int_0^t1_{\{x_{2r-2}^{m+1}< X_{s-}\leq
x_{2r-1}^{m+1}\}}dB_s+(K_1(t,x_{2r-1}^{m+1})-K_1(t,x_{2r-2}^{m+1}))\\
&&\hskip 0.5cm+(K_2(t,x_{2r-1}^{m+1})-K_2(t,x_{2r-2}^{m+1}))+(K_3(t,x_{2r-1}^{m+1})-K_3(t,x_{2r-2}^{m+1}))\Big]\\
&&\cdot\Big[\varphi_t(x_{2l-1}^{m+1})-\varphi_t(x_{2l-2}^{m+1})-b\int_0^t1_{\{x_{2l-2}^{m+1}<
X_{s-}\leq x_{2l-1}^{m+1}\}}ds\\
&&\hskip 0.5cm -\sigma\int_0^t1_{\{x_{2l-2}^{m+1}< X_{s-}\leq
x_{2l-1}^{m+1}\}}dB_s+(K_1(t,x_{2l-1}^{m+1})-K_1(t,x_{2l-2}^{m+1}))\\
&&\hskip
0.5cm+(K_2(t,x_{2l-1}^{m+1})-K_2(t,x_{2l-2}^{m+1}))+(K_3(t,x_{2l-1}^{m+1})-K_3(t,x_{2l-2}^{m+1}))\Big].
\end{eqnarray*}
Firstly, from (\ref{hz10}), (\ref{levy11a}), (\ref{levy4}), the
Cauchy-Schwarz inequality and the quadratic variation of stochastic
integrals, we have
\begin{eqnarray}\label{chunrong21}
&&\Big|E\big(\varphi_t(x_{2r-1}^{m+1})-\varphi_t(x_{2r-2}^{m+1})-b\int_0^t1_{\{x_{2r-2}^{m+1}<
X_{s-}\leq x_{2r-1}^{m+1}\}}ds-\sigma\int_0^t1_{\{x_{2r-2}^{m+1}<
X_{s-}\leq x_{2r-1}^{m+1}\}}dB_s\big)\nonumber\\
&&\hskip 0.5cm\cdot
\big(\varphi_t(x_{2l-1}^{m+1})-\varphi_t(x_{2l-2}^{m+1})-b\int_0^t1_{\{x_{2l-2}^{m+1}<
X_{s-}\leq x_{2l-1}^{m+1}\}}ds-\sigma\int_0^t1_{\{x_{2l-2}^{m+1}<
X_{s-}\leq x_{2l-1}^{m+1}\}}dB_s\big)\Big|\nonumber\\
%&\leq&E\big|\varphi_t(x_{2r-1}^{m+1})-\varphi_t(x_{2r-2}^{m+1})\big|\cdot\big|\varphi_t(x_{2l-1}^{m+1})-\varphi_t(x_{2l-2}^{m+1})\big|\nonumber\\
%&&+E\big|\varphi_t(x_{2r-1}^{m+1})-\varphi_t(x_{2r-2}^{m+1})\big|\cdot\big|b\int_0^t1_{\{x_{2l-2}^{m+1}\leq
%X_s<x_{2l-1}^{m+1}\}}ds\big|\nonumber\\
%&&+
% E\big|\varphi_t(x_{2r-1}^{m+1})-\varphi_t(x_{2r-2}^{m+1})\big|\cdot\big|\sigma\int_0^t1_{\{x_{2l-2}^{m+1}\leq
%X_s<x_{2l-1}^{m+1}\}}dB_s\big|\nonumber\\
%&&+
%E\big|\varphi_t(x_{2l-1}^{m+1})-\varphi_t(x_{2l-2}^{m+1})\big|\cdot\big|b\int_0^t1_{\{x_{2r-2}^{m+1}\leq
%X_s<x_{2r-1}^{m+1}\}}ds\big|\nonumber\\
%&&+
%E\big|\varphi_t(x_{2l-1}^{m+1})-\varphi_t(x_{2l-2}^{m+1})\big|\cdot\big|\sigma\int_0^t1_{\{x_{2r-2}^{m+1}\leq
%X_s<x_{2r-1}^{m+1}\}}dB_s\big|\nonumber\\
%&&+b^2 E\big|\int_0^t1_{\{x_{2r-2}^{m+1}\leq
%X_s<x_{2r-1}^{m+1}\}}ds\int_0^t1_{\{x_{2l-2}^{m+1}\leq
%X_s<x_{2l-1}^{m+1}\}}ds\big|\nonumber\\
%&&+\sigma^2 E\big|\int_0^t1_{\{x_{2r-2}^{m+1}\leq
%X_s<x_{2r-1}^{m+1}\}}1_{\{x_{2l-2}^{m+1}\leq
%X_s<x_{2l-1}^{m+1}\}}ds\big|\nonumber\\
&\leq&
c(t)\Big[(1+2b+b^2)(x_{2r-1}^{m+1}-x_{2r-2}^{m+1})(x_{2l-1}^{m+1}-x_{2l-2}^{m+1})
+\sigma(x_{2r-1}^{m+1}-x_{2r-2}^{m+1})(x_{2l-1}^{m+1}-x_{2l-2}^{m+1})^{1\over
2}\nonumber\\
&&\hskip
0.5cm+\sigma(x_{2l-1}^{m+1}-x_{2l-2}^{m+1})(x_{2r-1}^{m+1}-x_{2r-2}^{m+1})^{1\over
2}\Big]\nonumber\\
&&+\sigma^2 E\big|\int_0^t1_{\{x_{2r-2}^{m+1}\leq
X_s<x_{2r-1}^{m+1}\}}1_{\{x_{2l-2}^{m+1}\leq
X_s<x_{2l-1}^{m+1}\}}ds\big|\nonumber\\
&\leq&C\Big[({1\over{2^{m+1}}})^2
w_1(x',x'')^2+({1\over{2^{m+1}}})^{3\over 2} w_1(x',x'')^{3\over
2}\Big]\nonumber\\
&&+\sigma^2E\big|\int_0^t1_{\{x_{2r-2}^{m+1}\leq
X_s<x_{2r-1}^{m+1}\}}1_{\{x_{2l-2}^{m+1}\leq
X_s<x_{2l-1}^{m+1}\}}ds\big|\nonumber\\
&\leq &\cases {C({1\over{2^{m+1}}})^{3\over 2}, {\rm \ if } \ r\neq
l,\cr C{1\over{2^{m+1}}},\ \  {\rm \ if}\ r=l.}
\end{eqnarray}
Here $C$ is a generic constant and also depends on $t$, $b$,
$\sigma$, $w_1(x',x'')$. Secondly, recall the fact that
$E(P_.M_.)=0$, if $P$ is a process of bounded variation and $M$ is a
martingale with mean $0$ and at least one of $M$ and $P$ is
continuous. Note here $K_1$ is a process of bounded variation.
Recall also that the cross-variation of $\int_0^t1_{\{a_{i}<
X_{s-}\leq a_{i+1}\}}dB_s$ and the jump parts such as
$\big(K_2(t,a_{j+1})-K_2(t,a_{j})\big)$ are zero. So we have
\begin{eqnarray*}
&&E\int_0^t1_{\{a_{i}< X_{s-}\leq a_{i+1}\}}dB_s
\cdot\big(K_1(t,a_{j+1})-K_1(t,a_{j})\big)=0,\\
&&E\int_0^{t}1_{\{a_{i}< X_{s-}\leq a_{i+1}\}}dB_s\cdot\big(K_2(t,a_{j+1})-K_2(t,a_{j})\big)=0, \\
&&E\int_0^{t}1_{\{a_{i}< X_{s-}\leq
a_{i+1}\}}ds\cdot\big(K_2(t,a_{j+1})-K_2(t,a_{j})\big)=0, \\
%&&E\int_0^{t}1_{\{x_{2l-2}^{m+1}< X_{s-}\leq
%x_{2l-1}^{m+1}\}}dB_s\cdot\int_0^{t+}\int_{|y|<1}1_{\{x_{2r-2}^{m+1}<
%X_{s-}\leq x_{2r-1}^{m+1}\}}y \tilde N_p(dsdy)=0, \\
&&E\int_0^t1_{\{a_{i}< X_{s-}\leq a_{i+1}\}}dB_s
\cdot\big(K_3(t,a_{j+1})-K_3(t,a_j)\big)=0,\\
&&E\big(K_2(t,a_{i+1})-K_2(t,a_{i})\big)\cdot\big(K_3(t,a_{j+1})-K_3(t,a_{j})\big)=0.
\end{eqnarray*}
Thirdly, by Lemma \ref{lem8}, we can see that when $0\leq
\xi\leq{1\over 6}$,
\begin{eqnarray*}
\Big|E\big(K_2(t,a_{i+1})-K_2(t,a_{i})\big)\big(K_2(t,a_{j+1})-K_2(t,a_{j})\big)\Big|
\leq c(t,\sigma)[(a_{i+1}-a_i)^{1+2\xi}+(a_{j+1}-a_j)^{1+2\xi}].
%&&\Big|E\big(K_2(t,a_{i+1})-K_2(t,a_{i})\big)\int_0^{t+}\int_{|y|<1}
%1_{\{a_{j}< X_{s-}\leq a_{j+1}\}} y\tilde N_p(dsdy)\Big| \leq
%c(t,\sigma)[(a_{i+1}-a_i)^{4\over 3}+(a_{j+1}-a_j)^{4\over 3}].
\end{eqnarray*}
 For other terms,
%except for the last one in (\ref{1025}),
 by the Cauchy-Schwarz
inequality, it is easy to see that
\begin{eqnarray*}
&&\Big|E\big(\varphi_t(a_{i+1})-\varphi_t(a_{i})\big)\big(K_1(t,a_{j+1})-K_1(t,a_{j})\big)\Big|
\leq
c(t)[(a_{i+1}-a_i)(a_{j+1}-a_j)];\\
&&\Big|E\big(\varphi_t(a_{i+1})-\varphi_t(a_{i})\big)\big(K_2(t,a_{j+1})-K_2(t,a_{j})\big)\Big|
\leq c(t,\sigma)[(a_{i+1}-a_i)(a_{j+1}-a_j)^{1\over 2}];\\
&&\Big|E\big(\varphi_t(a_{i+1})-\varphi_t(a_{i})\big)\big(K_3(t,a_{j+1})-K_3(t,a_{j})\big)\Big|
\leq
c(t)[(a_{i+1}-a_i)(a_{j+1}-a_j)^{{1\over 2}+\xi};\\
&&\Big|E\int_0^t1_{\{a_i<X_{s-}\leq
a_{i+1}\}}ds\big(K_1(t,a_{j+1})-K_1(t,a_{j})\big)\Big|
\leq c(t,\sigma)[(a_{i+1}-a_i)(a_{j+1}-a_j)];\\
&&\Big|E\int_0^t1_{\{a_i<X_{s-}\leq
a_{i+1}\}}ds\big(K_3(t,a_{j+1})-K_3(t,a_{j})\big)\Big|
\leq c(t,\sigma)[(a_{i+1}-a_i)(a_{j+1}-a_j)^{{1\over 2}+\xi}];\\
&&\Big|E\big(K_1(t,a_{i+1})-K_1(t,a_{i})\big)\big(K_1(t,a_{j+1})-K_1(t,a_{j})\big)\Big|
\leq
c(t,\sigma)[(a_{i+1}-a_i)(a_{j+1}-a_j)];\\
&&\Big|E\big(K_1(t,a_{i+1})-K_1(t,a_{i})\big)\big(K_2(t,a_{j+1})-K_2(t,a_{j})\big)\Big|
\leq c(t,\sigma)[(a_{i+1}-a_i)(a_{j+1}-a_j)^{1\over 2}];\\
&&\Big|E\big(K_1(t,a_{i+1})-K_1(t,a_{i})\big)\big(K_3(t,a_{j+1})-K_3(t,a_{j})\big)\Big|
\leq
c(t,\sigma)[(a_{i+1}-a_i)(a_{j+1}-a_j)^{{1\over 2}+\xi}];\\
&&\Big|E\big(K_3(t,a_{i+1})-K_3(t,a_{i})\big)\big(K_3(t,a_{j+1})-K_3(t,a_{j})\big)\Big|
\leq c(t,\sigma)[(a_{i+1}-a_i)^{{1\over
2}+\xi}(a_{j+1}-a_j)^{{1\over 2}+\xi}].
\end{eqnarray*}
 Thus
\begin{eqnarray}\label{chunrong21}
\left|E\Big[\Delta_{2r-1}^{m+1}L_t^x\Delta_{2l-1}^{m+1}L_t^x\Big]\right|
\leq \cases {C({1\over{2^{m+1}}})^{1+2\xi}, {\rm \ if } \ r\neq
l,\cr C{1\over{2^{m+1}}},\ \  {\rm \ if}\ r=l.}
\end{eqnarray}
The other terms in (\ref {control7}) can be treated similarly.
Therefore
\begin{eqnarray*}
E|{\bf Z}(m+1)_{x_{k-1}^n,x_k^n}^2-{\bf Z}(m)_{x_{k-1}^n,x_k^n}^2|^2
\leq
C\Big[2^{m-n}({1\over{2^{m+1}}})^{1+2h}+2^{2(m-n)}({1\over{2^{m+1}}})^{{1+2\xi}+2h}\Big].
\end{eqnarray*}
Hence, for $2\leq \theta<3$, by Jensen's inequality,
\begin{eqnarray*}
E|{\bf Z}(m+1)_{x_{k-1}^n,x_k^n}^2-{\bf
Z}(m)_{x_{k-1}^n,x_k^n}^2|^{\theta\over 2} &\leq& \left(E|{\bf
Z}(m+1)_{x_{k-1}^n,x_k^n}^2-{\bf
Z}(m)_{x_{k-1}^n,x_k^n}^2|^2\right)^{\theta\over
4}\\
&\leq&
C\Big[2^{m-n}({1\over{2^{m+1}}})^{1+2h}+2^{2(m-n)}({1\over{2^{m+1}}})^{{1+2\xi}+2h}\Big]^{\theta\over
4}\\
&\leq& C\Big[2^{(m-n){\theta\over
4}}({1\over{2^{m+1}}})^{{\theta\over 4}+{1\over
2}h\theta}+2^{(m-n){\theta\over
2}}({1\over{2^{m+1}}})^{{{1+2\xi}\over 4}\theta+{1\over 2}h\theta}\Big]\\
&\leq& C\Big[({1\over{2^n}})^{\theta\over 4}({1\over{2^m}})^{{1\over
2}h\theta}+({1\over{2^n}})^{\theta\over 2}({1\over{2^m}})^{{1\over
2}h\theta-{{1-2\xi}\over 4}\theta}\Big],
\end{eqnarray*}
where $C$ is a generic constant and also depends on $\theta$, $h$,
$w_1(x',x'')$, and $c$. Note $\xi={{q-2}\over {2q}}+(3-q)\varepsilon$, so we get (\ref{may13}).  $\hfill\diamond$
 \vskip5pt
\begin{cor}\label{cornew1}
Under the same assumption as in Proposition \ref {propo2}, we have
\begin{eqnarray*}
\sup_m\sum_{n=1}^\infty n^{\gamma}\sum_{k=1}^{2^n}|{\bf
Z}(m)_{x_{k-1}^n,x_k^n}^2|^{\theta\over 2}<\infty\ \ \ a.s.
\end{eqnarray*}
\end{cor}
{\bf Proof:} From the Minkowski inequality,
\begin{eqnarray}\label{may16}
&&\left (\sum_{k=1}^{2^n}|{\bf Z}(m)_{x_{k-1}^n,x_k^n}^2|^{\theta\over 2}\right)^{2\over \theta}\nonumber\\
&\leq & \bigg(\sum_{k=1}^{2^n}|{\bf Z}(m)_{x_{k-1}^n,x_k^n}^2-{\bf Z}(m-1)_{x_{k-1}^n,x_k^n}^2|^{\theta\over 2}\bigg)^{2\over \theta}
%\nonumber\\
%&&
+
\bigg(\sum_{k=1}^{2^n}|{\bf Z}(m-1)_{x_{k-1}^n,x_k^n}^2-{\bf Z}(m-2)_{x_{k-1}^n,x_k^n}^2|^{\theta\over 2}\bigg)^{2\over \theta}\nonumber\\
&&+\cdots+\bigg(\sum_{k=1}^{2^n}|{\bf Z}(1)_{x_{k-1}^n,x_k^n}^2-{\bf
Z}(0)_{x_{k-1}^n,x_k^n}^2|^{\theta\over 2}\bigg)^{2\over \theta}
+\bigg(\sum_{k=1}^{2^n}|{\bf Z}(0)_{x_{k-1}^n,x_k^n}^2|^{\theta\over 2}\bigg)^{2\over \theta}\nonumber\\
&\leq & \sum_{m=1}^\infty \bigg(\sum_{k=1}^{2^n}|{\bf
Z}(m)_{x_{k-1}^n,x_k^n}^2-{\bf
Z}(m-1)_{x_{k-1}^n,x_k^n}^2|^{\theta\over 2}\bigg)^{2\over
\theta}+\bigg(\sum_{k=1}^{2^n}|{\bf
Z}(0)_{x_{k-1}^n,x_k^n}^2|^{\theta\over 2}\bigg)^{2\over \theta}.
\end{eqnarray}
Then it is easy to see from (\ref{may16}), Jensen's inequality,
(\ref{may12}), (\ref{may13}) and (\ref{cfeng3})
\begin{eqnarray*}
&&E\sup_m \sum_{n=1}^\infty n^{{2\over \theta}\gamma}\left(\sum_{k=1}^{2^n}|{\bf Z}(m)_{x_{k-1}^n,x_k^n}^2|^{\theta\over 2}\right)^{2\over \theta}\\
&\leq & E\sum_{m=1}^\infty \left [\sum_{n=1}^\infty n^{{2\over
\theta}\gamma}\left(\sum_{k=1}^{2^n}|{\bf
Z}(m)_{x_{k-1}^n,x_k^n}^2-{\bf
Z}(m-1)_{x_{k-1}^n,x_k^n}^2|^{\theta\over 2}\right)^{2\over \theta}
\right .\\
&& \left . \hskip3cm +\sum_{n=1}^\infty n^{{2\over
\theta}\gamma}E\bigg(\sum_{k=1}^{2^n}|{\bf
Z}(0)_{x_{k-1}^n,x_k^n}^2|^{\theta\over
2}\bigg)^{2\over \theta}\right ]\\
&\leq &\sum_{m=1}^\infty \left [\sum_{n=1}^\infty n^{{2\over
\theta}\gamma}\left(E\sum_{k=1}^{2^n}|{\bf
Z}(m)_{x_{k-1}^n,x_k^n}^2-{\bf
Z}(m-1)_{x_{k-1}^n,x_k^n}^2|^{\theta\over 2}\right)^{2\over
\theta}\right .\\
&& \hskip3cm
\left .+\sum_{n=1}^\infty n^{{2\over \theta}\gamma}(2^{n+{{(-1-2n)\theta}\over 2}}w(x',x'')^{h\theta})^{2\over \theta}\right ]\\
&\leq & C\sum_{n=m}^{\infty} n^{{2\over
\theta}\gamma}\sum_{m=1}^{\infty}
({1\over {2^{n+m}}})^{({{h\theta-1}\over 2}){2\over \theta}}\\
&&+C\sum_{n=1}^{\infty}n^{{2\over \theta}\gamma}\sum_{m=1}^\infty
\Big[({1\over{2^n}})^{({\theta\over 4}-{1\over 2}){2\over
\theta}}({1\over{2^{m}}})^{({{1\over 2}h\theta}-{1\over 2}){2\over
\theta}}+({1\over{2^n}})^{1-{2\over
\theta}}({1\over{2^{m}}})^{({{{3-q}\over 2}{\varepsilon\theta}}){2\over \theta}}\Big]\\
&&+C\sum_{n=1}^{\infty }n^{{2\over \theta} \gamma}({1\over
{2^n}})^{{2(\theta-1)}\over
\theta}\\
 &< &\infty,
\end{eqnarray*}
as $2<\theta<3$, $h\theta>1$, where $C$ is a generic constant and
also depends on $\theta$, $h$, $w_1(x',x'')$, and $c$. Therefore,
\begin{eqnarray*}
\sup_m \sum_{n=1}^\infty n^{{\theta\over
2}\gamma}\left(\sum_{k=1}^{2^n}|{\bf
Z}(m)_{x_{k-1}^n,x_k^n}^2|^{\theta\over 2}\right)^{2\over
\theta}<\infty\ \ a.s.
\end{eqnarray*}
 However, it is easy to see as $\theta>2$,
\begin{eqnarray*}
&&\left(\sup_m \sum_{n=1}^\infty n^{\gamma}\sum_{k=1}^{2^n}|{\bf
Z}(m)_{x_{k-1}^n,x_k^n}^2|^{\theta\over 2}\right)^{2\over\theta}\leq
\sup_m \sum_{n=1}^\infty n^{{2\over
\theta}\gamma}\left(\sum_{k=1}^{2^n}|{\bf
Z}(m)_{x_{k-1}^n,x_k^n}^2|^{\theta\over 2}\right)^{2\over
\theta}<\infty\ \ a.s.
\end{eqnarray*}
So the claim follows. $\hfill\diamond$
 \vskip5pt
\begin{thm}\label{theor2}
Let $L_t^x$ be the local time of the time homogeneous L\'evy process
$X_t$ given by (\ref{tanakaf1}). Assume $g$ is a continuous function of finite 
$q$-variation with a real number $2\leq q<3$,
and the L\'evy measure $n(dy)$ satisfies (\ref{levynew36}), $\sigma\neq 0$. Then for
any $\theta\in (q,3)$, the continuous process $Z_x=(L_t^x, g(x))$
satisfying (\ref{control5}), there exists a unique ${\bf Z}^i$ on
$\Delta$ taking values in $(R^2)^{\otimes i}$ ($i=1,2$) such that
\begin{eqnarray*}
\sum_{i=1}^2\sup_D\bigg(\sum_l|{\bf Z}(m)_{x_{l-1},x_l}^i-{\bf
Z}_{x_{l-1},x_l}^i|^{\theta\over i}\bigg)^{i\over \theta}\to 0 ,
\end{eqnarray*}
both almost surely and in $L^1(\Omega, \cal F, P)$ as $m\to\infty$.
In particular, ${\bf Z}=(1,{\bf Z}^1, {\bf Z}^2)$ is the canonical
geometric rough path associated to $Z_{.}$, and ${\bf
Z}^1_{a,b}=Z_b-Z_a$.
\end{thm}
{\bf Proof:} The convergence of ${\bf Z}(m)^1$ to ${\bf Z}^1$ is
actually the result of Theorem \ref{theor1}. In the following
 we will prove ${\bf Z}(m)_{a,b}^2$ converges in the $\theta$-variation
distance. By Proposition 4.1.2 in \cite{terry},
\begin{eqnarray*}
&&E\sup_D\sum_l|{\bf Z}(m+1)_{x_{l-1},x_l}^2-{\bf
Z}(m)_{x_{l-1},x_l}^2|^{\theta\over
2}\\
&\leq& C(\theta, \gamma)E\left(\sum_{n=1}^\infty
n^{\gamma}\sum_{k=1}^{2^n}|{\bf Z}(m+1)_{x_{k-1}^n,x_k^n}^1-{\bf
Z}(m)_{x_{k-1}^n,x_k^n}^1|^{\theta}\right)^{1\over
2}\\
&&\hskip 1.5cm\cdot\left(\sum_{n=1}^\infty
n^{\gamma}\sum_{k=1}^{2^n}\Big(|{\bf
Z}(m+1)_{x_{k-1}^n,x_k^n}^1|^{\theta}+|{\bf
Z}(m)_{x_{k-1}^n,x_k^n}^1|^{\theta}\Big)\right)^{1\over
2}\\
&&+C(\theta, \gamma)E\sum_{n=1}^\infty
n^{\gamma}\sum_{k=1}^{2^n}|{\bf Z}(m+1)_{x_{k-1}^n,x_k^n}^2-{\bf
Z}(m)_{x_{k-1}^n,x_k^n}^2|^{\theta\over
2}\\
&:= &A+B.
\end{eqnarray*}
We will  estimate part A, B respectively. First from the
Cauchy-Schwarz inequality, (\ref{newww1}) and (\ref{cfeng2}), we
know
\begin{eqnarray*}
A&\leq& C\bigg(E\sum_{n=1}^\infty n^{\gamma}\sum_{k=1}^{2^n}\bigg
(|{\bf Z}(m+1)_{x_{k-1}^n,x_k^n}^1-{\bf
Z}_{x_{k-1}^n,x_k^n}^1|^{\theta} +|{\bf
Z}(m)_{x_{k-1}^n,x_k^n}^1-{\bf Z}_{x_{k-1}^n,x_k^n}^1|^{\theta}\bigg
)\bigg)^{1\over
2}\\
&&\hskip 1cm\cdot\bigg(E\sum_{n=1}^\infty
n^{\gamma}\sum_{k=1}^{2^n}\bigg (|{\bf
Z}(m+1)_{x_{k-1}^n,x_k^n}^1|^{\theta}+|{\bf
Z}(m)_{x_{k-1}^n,x_k^n}^1|^{\theta}\bigg )\bigg)^{1\over
2}\\
&\leq& C({1\over {2^m}})^{{h\theta-1}\over 4}\bigg(\sum_{n=1}^\infty
n^\gamma({1\over {2^n}})^{{h\theta-1}}\bigg )^{1\over 2}.
\end{eqnarray*}
 Secondly from
Proposition \ref{propo1} and Proposition \ref{propo2}, we know
\begin{eqnarray*}
B&\leq& C\sum_{n=m}^\infty n^\gamma({1\over
{2^{m+n}}})^{{h\theta-1}\over
2}+C\bigg[\sum_{n=1}^{m-1}n^\gamma({1\over {2^{n}}})^{{\theta\over
4}-1}({1\over {2^{m}}})^{{{h\theta}\over
2}}+\sum_{n=1}^{m-1}n^\gamma({1\over {2^{n}}})^{{\theta\over
2}-1}({1\over {2^{m}}})^{{{{3-q}\over 2}{\varepsilon\theta}}}\bigg]\\
&\leq& C\bigg[({1\over {2^m}})^{{h\theta-1}\over 2}+({1\over
{2^m}})^{{h\theta-1}\over 2}+({1\over {2^{m}}})^{{{{3-q}\over
2}{\varepsilon\theta}}}\bigg],
\end{eqnarray*}
as $q<\theta<3$, and $h\theta>1$. So
\begin{eqnarray*}
E\sup_D\sum_l|{\bf Z}(m+1)_{x_{l-1},x_l}^2-{\bf
Z}(m)_{x_{l-1},x_l}^2|^{\theta\over 2}\leq C[({1\over
{2^m}})^{{h\theta-1}\over 4}+({1\over {2^{m}}})^{{{{3-q}\over
2}{\varepsilon\theta}}}].
\end{eqnarray*}
Similar to the proof of Theorem \ref{theor1}, we can easily deduce
that $({\bf Z}(m)^2)_{m\in N}$ is a Cauchy sequence in the
$\theta$-variation distance. So when $m\to\infty$, it has a limit,
denote it by ${\bf Z}^2$. And from the completeness under the
$\theta$-variation distance (Lemma 3.3.3 in \cite{terry}), ${\bf
Z}^2$ is also of finite $\theta$-variation. The theorem is proved.
$\hfill\diamond$
 \vskip5pt

\begin{rmk}
We would like to point out that the above method does not seem to
work for two arbitrary functions $f$ of $p$-variation and $g$ of
$q$-variation ($2<p,q<3$) to define a rough path $Z_x=(f(x),g(x))$.
However the special property (\ref{chunrong21}) of local times makes
our analysis work. A similar method was used in \cite{terry} for
fractional Brownian motion with the help of long-time memory. Here
(\ref{chunrong21}) serves a similar role of the long-time memory as
in \cite{terry}.
\end{rmk}

As local time $L_t^x$ has a compact support in $x$ for each $\omega$
and $t$, so we can define integral of local time directly in $R$.
For this, we take $[x^{\prime},x^{\prime\prime}]$ covering the
support of $L_t^x$. From Chen's identity, it's easy to know that for
any $(a,b)\in \Delta$,
\begin{eqnarray*}
{\bf Z}_{a,b}^2=\lim_{m(D_{[a,b]})\to 0}\sum_{i=0}^{r-1}({\bf
Z}_{x_i,x_{i+1}}^2+ {\bf Z}_{a,x_i}^1\otimes {\bf
Z}_{x_i,x_{i+1}}^1).
\end{eqnarray*}
In particular,
\begin{eqnarray*}
%\int_{a<x_1<x_2<b}d g(x_1)dL_t(x_2)&=&
({\bf Z}_{a,b}^2)_{2,1} &=& \lim_{m(D_{[a,b]})\to
0}\sum_{i=0}^{r-1}(({\bf Z}_{x_i,x_{i+1}}^2)_{2,1}+
({\bf Z}_{a,x_i}^1\otimes {\bf Z}_{x_i,x_{i+1}}^1)_{2,1})\\
&=&\lim_{m(D_{[a,b]})\to 0}\sum_{i=0}^{r-1}(({\bf
Z}_{x_i,x_{i+1}}^2)_{2,1}+(g(x_i)-g(a))(L_t^{x_{i+1}}- L_t^{x_i}))
\end{eqnarray*}
exists. Here $({\bf Z}_{x_i,x_{i+1}}^2)_{2,1}$ means lower-left element
of the $2\times 2$ matrix ${\bf
Z}_{x_i,x_{i+1}}^2$.
It turns out that
\begin{eqnarray*}
&&\lim_{m(D_{[a,b]})\to 0}\sum_{i=0}^{r-1}(({\bf
Z}_{x_i,x_{i+1}}^2)_{2,1}+g(x_i)(
L_t^{x_{i+1}}-L_t^{x_i}))\\
&=&\lim_{m(D_{[a,b]})\to 0}\sum_{i=0}^{r-1}(({\bf
Z}_{x_i,x_{i+1}}^2)_{2,1}+(g(x_i)-g(a))(L_t^{x_{i+1}}-
L_t^{x_i}))+g(a)(L_t^{b}-L_t^{a})\nonumber
\end{eqnarray*}
exists. Denote this limit by $\int_a^b g(x)dL_t^x$. Similarly, we
can define $\int_a^b L_t^xdL_t^x$. To verify the latter integral is
well defined, we only need to consider the case $q=2$. Then it is
easy to see under condition (\ref{levynew37}), $\int_a^b
L_t^xdL_t^x$ is defined as a rough path integral. Therefore we have
the following corollary.
\begin{cor}
Assume all conditions of Theorem \ref{theor2}, but the L\'evy
measure satisfies (\ref{levynew37}). Then the local time $L_t^x$ is
a geometrical rough path of roughness $p$ in $x$ for any $t\geq 0$
a.s. for any $p>2$, and $(a,b)\in \Delta$,
\begin{eqnarray*}
\int_a^b L_t^xdL_t^x=\lim_{m(D_{[a,b]})\to 0}\sum_{i=0}^{r-1}(({\bf
Z}_{x_i,x_{i+1}}^2)_{1,1}+L(x_i)(L_t^{x_{i+1}}-L_t^{x_i})).
\end{eqnarray*}
Moreover, if $g$ is a continuous function with bounded
$q$-variation, $2\leq q<3$, and the L\'evy measure satisfies
(\ref{levynew36}), then the integral $\int_a^b g(x)d L_t^x$ is
defined by
\begin{eqnarray}\label{cfeng12}
\int_a^b g(x)dL _t(x)=\lim_{m(D_{[a,b]})\to 0}\sum_{i=0}^{r-1}(({\bf
Z}_{x_i,x_{i+1}}^2)_{2,1}+g(x_i)(L_t^{x_{i+1}}-L_t^{x_i})).
\end{eqnarray}
\end{cor}

% Note it is clear to us that the
%Riemann sum $\sum\limits _{i=0}^{r-1}g(x_i)( L_t^{x_{i+1}}-L_t^{x_i})$
%and $\sum\limits _{i=0}^{r-1} L(x_i)(L_t^{x_{i+1}}-L_t^{x_i})$
%themselves do not have limits as $m(D_{[a,b]}) \to 0$. But the rough
%path integration theory provides a way to define such a path
%integral.

\section{Continuity of the rough path integrals and applications to extensions of It$\hat{\rm o}$'s formula}

In this section we will apply the Young integral and rough path
integral of local time defined in sections 2 and 3 to prove a useful
extension to It$\hat{\rm o}$'s formula. First we consider some
convergence result of the rough path integrals.

  Let
 $Z_j(x):=(L_t^x,g_j(x))$, where
 $g_j(\cdot)$ is of bounded $q$-variation uniformly in $j$ for $2\leq
 q<3$, and when $j\to\infty$, $g_j(x)\to g(x)$ for all $x\in R$. Assume the L\'evy
measure satisfies (\ref{levynew36}). Repeating the above
 argument, for each $j$, we can find the canonical geometric rough path ${\bf Z}_j=(1, {\bf Z}_j^1, {\bf Z}_j^2)$ associated to
 $Z_j$, and the smooth rough path ${\bf Z}_j(m)=(1, {\bf Z}_j(m)^1,
 {\bf Z}_j(m)^2)$. Actually,
 % by the result about ${\bf Z}^1$ in \cite{terry}, it is easy to see
  %that $({\bf Z}_j)_{a,b}^1=(L_t^b-L_t^a,
 %g_j(b)-g_j(a))$, ${\bf Z}_{a,b}^1=(L_t^b-L_t^a,g(b)-g(a))$,
 $({\bf Z}_j)_{a,b}^1\to {\bf Z}_{a,b}^1$ in the
sense
 of the
 uniform topology, and also in the sense of the $\theta$-variation
 topology. As for $({\bf Z}_j)_{a,b}^2$, we can easily see that
 \begin{eqnarray}\label{may2}
d_{2,\theta}(({\bf Z}_j)^2, {\bf Z}^2) \leq d_{2,\theta}(({\bf
Z}_j)^2,({\bf Z}_j(m))^2)+d_{2,\theta}(({\bf Z}_j(m))^2,{\bf
Z}(m)^2)+d_{2,\theta}({\bf Z}(m)^2,{\bf Z}^2).
\end{eqnarray}
From Theorem \ref{theor2}, we know that $d_{2,\theta}({\bf
Z}(m)^2,{\bf Z}^2)\to 0$ as $m\to \infty$. Moreover, it is not
difficult to see from the proofs of Propositions \ref{propo1},
\ref{propo2}, and Theorem \ref{theor2}, $d_{2,\theta}(({\bf
Z}_j)^2,({\bf Z}_j(m))^2)\to 0$ as $m\to \infty$ uniformly in $j$.
So for any given $\varepsilon>0$, there exists an $m_0$ such that
when $m\geq m_0$, $d_{2,\theta}({\bf Z}(m)^2,{\bf
Z}^2)<{\varepsilon\over 3}$, $d_{2,\theta}(({\bf Z}_j)^2,({\bf
Z}_j(m))^2)<{\varepsilon\over 3}$ for all $j$. In particular,
$d_{2,\theta}({\bf Z}(m_0)^2,{\bf Z}^2)<{\varepsilon\over 3}$,
$d_{2,\theta}(({\bf Z}_j)^2,({\bf Z}_j(m_0))^2)<{\varepsilon\over
3}$ for all $j$. It's easy to prove for such $m_0$,
$d_{2,\theta}(({\bf Z}_j(m_0))^2,{\bf Z}(m_0)^2)< {\varepsilon\over
3}$ for sufficiently large $j$. Replacing $m$ by $m_0$ in
(\ref{may2}), we can get $d_{2,\theta}(({\bf Z}_j)^2, {\bf
Z}^2)<\varepsilon$ for sufficiently large $j$. Then  by
(\ref{cfeng12}) and the definition of $\int_a^b g_j(x)d L_t^x$, we
know that $\int_a^b g_j(x)dL_t^x\to \int_a^b g(x)d L_t^x$ as $j\to
\infty$. Similarly, we can see from the last section,  when we consider ${\bf Z}_t(m)=(1, {\bf Z}_t(m)^1,
 {\bf Z}_t(m)^2)$, $d_{2,\theta}(({\bf
Z}_t)^2,({\bf Z}_t(m))^2)\to 0$, as $m\to\infty$ uniformly in $t\in [0,T]$, for any $T>0$. Therefore we can also conclude that ${\bf Z}^2_t$ is continuous in $t$ in the $d_{2,\theta}$ topology.  Note now that the local time $L_t^x$ has a compact support
in $x$ a.s.  So it is easy to see from taking
$[x^{\prime},x^{\prime\prime}]$ covering the support of $L_t^x$ that
the above construction of the integrals and the convergence can work
for the integrals on $R$.
 Therefore we have
\begin{prop}\label{may5}
 Let
 $Z_j(x):=(L_t^x,g_j(x))$, $Z(x):=(L_t^x,g(x))$, where
 $g_j(\cdot)$, $g(\cdot)$ are continuous and of bounded $q$-variation uniformly in $j$, $2\leq
 q<3$,  and the L\'evy
measure $n(dy)$ satisfies (\ref{levynew36}), $\sigma \ne 0$. Assume  $g_j(x)\to
g(x)$ as $j\to\infty$ for all $x\in R$. Then as $j\to \infty$,
 ${\bf Z}_j(\cdot)\to {\bf Z}(\cdot)$ a.s. in the $\theta$-variation distance.
In particular, as $j\to \infty$, $\int_{-\infty}^\infty
g_j(x)dL_t^x \to \int_{-\infty}^\infty
  g(x)dL_t^x a.s.$ Similarly, $Z_t(\cdot)$ is continuous in $t$ in the $\theta$-topology. In particular,  $ \int_{-\infty}^\infty
  g(x)dL_t^x$ is a continuous function of $t$ a.s. 
\end{prop}

Now for any $g$ being continuous and of bounded $q$-variation
($2\leq q <3$), define
\begin{eqnarray}\label{levynew62}
g_j(x)=\int_{-\infty}^\infty k^j(x-y)g(y)dy,
\end{eqnarray}
where $k^j$ is the mollifier given by
\begin{eqnarray*}\label{smooth}
k^j(x)=\cases {cj{\rm e}^{{1\over (jx-1)^2-1}}, {\rm \  if } \ x\in
(0,{2\over j}),\cr 0, \ \ \ \ \ \ \ \ \ \ \ \ \  {\rm otherwise.}}
\end{eqnarray*}
Here $c$ is a constant such that $\int _0^2k^j(x)dx=1$. It is well
known that $g_j$ is a smooth function and $g_j(x)\to g(x)$ as $j\to
\infty$ for each $x$. So the integral $\int_{-\infty}^\infty g_j(x)
dL_t^x$ is a Riemann integral for the smooth function $g_j(x)$.
Moreover, Proposition \ref{may5} guarantees that
$\int_{-\infty}^\infty g_j(x)dL_t^x \to \int_{-\infty}^\infty
  g(x)dL_t^x$ a.s.

 In the following, we will show that Proposition \ref{may5} is true for $g$ being of bounded $q$-variation ($2\leq q<3$)
 without assuming $g$ being continuous. Note that a function with bounded $q$-variation
  ($q\geq 1$) may have at most countable discontinuities. Using the
  method in \cite{williams}, we will define the rough path integral $\int_{x'}^{x''}
  L_t^{x}dg(x)$. Here we assume $g(x)$ is c$\grave{a}$dl$\grave{a}$g in $x$.

  First we can define a map
  $$\tau_\delta(\cdot): [x',x'']\to [x', x''+\delta\sum\limits_{n=1}^\infty |j(x_n)|^q],$$
  in the following way:
  $$\tau_\delta(x)=x+\delta\sum\limits_{n=1}^\infty |j(x_n)|^q 1_{\{x_n\leq x\}}(x),$$
  where $j(x_i):=G(x_i)-G(x_i-)$, $\{x_i\}_{i=1}^\infty$ are the discontinuous points of $G$ inside $[x',x'']$, $\delta>0$.
  The map $\tau_\delta(\cdot): [x',x'']\to [x', \tau_\delta(x'')]$ extends
  the space interval into the one where we can define the continuous path
  $G_\delta(y)$ from a c$\grave{a}$dl$\grave{a}$g path $G$ by:
  \begin{eqnarray}\label{levy52}
G_\delta(y)=
  \left\{
\begin{array}{ll}
G(x) &{\rm if}  \ y=\tau_\delta (x),\\
G(x_n-)+(y-\tau_\delta(x_n-))j(x_n)\delta^{-1}|j(x_n)|^{-q} &{\rm
if} \ y\in [\tau_\delta(x_n-),\tau_\delta(x_n)).
\end{array}
\right.
\end{eqnarray}
Take $G$ to be $g$ and $L_t$, we can define $g_\delta$ and
$L_{t,\delta}$ respectively. As $L_t^x$ is continuous, we can easily
see that  $ L_{t,\delta}(y):={L}_{t,\delta}(\tau_\delta(x))=L_t^x$.
\begin{thm}
Let $g(x)$ be a c$\grave{a}$dl$\grave{a}$g path with bounded
$q$-variation ($2\leq q<3$), and the L\'evy
measure $n(dy)$ satisfies (\ref{levynew36}), $\sigma \ne 0$. Then
\begin{eqnarray}\label{jump1}
\int_{x'}^{x''} L_t^{x}dg(x)=\int_{x'}^{\tau_\delta(x'')}
{L}_{t,\delta}(y)dg_\delta(y).
\end{eqnarray}
\end{thm}
{\bf Proof:}  First it is easy to see that the integral $\int_{x'}^{\tau_\delta(x'')} {L}_{t,\delta}(y)dg_\delta(y)$
is a rough path integral that can be defined by the method of last section.
Now note that at any discontinuous point $x_r$,
\begin{eqnarray*}
\int_{x_r-}^{x_r}  L_t^{x}dg(x)=L_t(x_r) (g(x_r)-g(x_r-))
\end{eqnarray*}
and
\begin{eqnarray*}
\sum_{r}((Z_\delta)^2_{\tau_\delta(x_r-),\tau_\delta(x_r)})_{2,1}
&=&\sum_r\int_{\tau_\delta(x_r-)}^{\tau_\delta(x_r)} (
{L}_{t,\delta}(y)-
{L}_{t,\delta}(\tau_\delta(x_r-)))dg_\delta(y)=0,
\end{eqnarray*}
where $Z_\delta(y):=( {L}_{t,\delta}(y), g_\delta(y))$. Thus
\begin{eqnarray*}
&&\sum_{r}L_t^{x_r}(g(x_r)-g(x_r-))\\
&=&\sum_{r}
{L}_{t,\delta}(\tau_\delta(x_r-))(g_\delta(\tau_\delta(x_r))-g_\delta(\tau_\delta(x_r-)))\\
&=&\sum_{r}\left[
{L}_{t,\delta}(\tau_\delta(x_r-))(g_\delta(\tau_\delta(x_r))-g_\delta(\tau_\delta(x_r-)))+((Z_\delta)^2_{\tau_\delta(x_r-),\tau_\delta(x_r)})_{2,1}\right]\\
&<&\infty,
\end{eqnarray*}
so
\begin{eqnarray*}
\int_{\tau_\delta(x_r-)}^{\tau_\delta(x_r)}
{L}_{t,\delta}(y)dg_\delta(y)&=&
{L}_{t,\delta}(\tau_\delta(x_r-))(g_\delta(\tau_\delta(x_r))-g_\delta(\tau_\delta(x_r-)))\\
&=& {L}_t(x_r)(g(x_r)-g(x_{r-})).
\end{eqnarray*}
Thus
\begin{eqnarray*}
\int_{x_r-}^{x_r}
L_t^{x}dg(x)=\int_{\tau_\delta(x_r-)}^{\tau_\delta(x_r)}
{L}_{t,\delta}(y)dg_\delta(y).
\end{eqnarray*}
Now define $g(x)=\tilde g(x)+h(x)$, where $h(x)=\sum\limits_{x_r\leq
x}(g(x_r)-g(x_{r-}))$. Then $\tilde g$ is the continuous part of $g$
and $h$ is the jump part of $g$. Moreover, $\tilde g$ satisfies the
$q$-variation condition. So $\int_{x'}^{x''} L_t (x) d\tilde g(x)$
can be well defined as in the last section. For $h$, we can define
$h_{\delta}$ by taking $G=h$ in (\ref{levy52}). So the integral
$\int_{x'}^{x"}L_t^x  d h(x)$ can be well defined by the
followings:
\begin{eqnarray*}
\int_{x'}^{\tau_\delta(x'')} {L}_{t,\delta}(y)dh_\delta(y)
&=&\sum_r\int_{\tau_\delta(x_{r-})}^{\tau_\delta(x_{r})}
{L}_{t,\delta}(y)dh_\delta(y)
=\sum_r L_t(x_r)(h(x_r)-h(x_{r-}))\\
&=&\sum_r L_t(x_r)(g(x_r)-g(x_{r-}))=\sum_r\int_{x_{r-}}^{x_{r}}
L_t^{x}dh(x)=\int_{x'}^{x''}L_t^{x}dh(x).
\end{eqnarray*}
Therefore
\begin{eqnarray*}
\int_{x'}^{x''}L_t^x dg(x)&=&\int_{x'}^{x''}L_t^x d\tilde
g(x)+\int_{x'}^{x''}L_t^x dh(x)\\
&=&\int_{x'}^{\tau_\delta(x'')} {L}_{t,\delta}(y)d\tilde
g_\delta(y)+\int_{x'}^{\tau_\delta(x'')}
{L}_{t,\delta}(y)dh_\delta(y)\\
&=&\int_{x'}^{\tau_\delta(x'')} {L}_{t,\delta}(y)dg_\delta(y).
\end{eqnarray*}
$\hfill\diamond$
% But
%\begin{eqnarray*}
%\int_{x'}^{\tau^\delta(x'')}
%{L}^\delta_t(y)dg^\delta(y)=\lim_{m(D)\to 0}\sum_i\left[
%{L}^\delta_t(y_{i-1})(g^{\delta}(y_i)-g^\delta
%(y_{i-1}))+((Z^\delta)^2_{y_{i-1},y_i})_{2,1}\right],
%\end{eqnarray*}
%is the rough path integral, . For the discontinuities $x_r,
%r=1,2,\cdots$, we have
%
%In other intervals excluding $x_r,r=1,2,\cdots$, (\ref{jump1}) holds
%obviously. If the discontinuities are dense, we just define the
%integral
%$$\int_{x'}^{x''} {L}_t(x)dg(x):=\int_{x'}^{\tau^\delta(x'')}
%{L}^\delta_t(y)dg^\delta(y). $$ So the claim follows.
 \vskip5pt
Similarly to Proposition \ref{may5}, we have
\begin{prop}
Under the condition of Proposition \ref{may5}, as $j\to \infty$,
$\int_{-\infty}^\infty L_t^x dg_j(x) \to \int_{-\infty}^\infty
 L_t^x dg(x)\ \  a.s.$ for such $g$ with bounded
 $q$-variation ($2\leq q <3$).
\end{prop}
{\bf Proof:} Define $F_j(x):=(g_j-g)(x)$, so $F_j(x)\to 0$ as $j\to
\infty$, for all $x$. It's easy to see that $F_{j,\delta}(x)\to 0$ as
$j\to \infty$, for all $x$. From the above theorem and Proposition
\ref{may5}, we have
\begin{eqnarray*}
\int_{x'}^{x''}L_t^x d(g_j-g)(x)=\int_{x'}^{x''}L_t^x
dF_j(x)=\int_{x'}^{\tau_\delta(x'')}L_{t,\delta}dF_{j,\delta}(y)\to
0,\ \ as\ j\to \infty.
\end{eqnarray*}
Then the proposition follows easily. $\hfill\diamond$
 \vskip5pt
%\begin{eqnarray*}
%\int_{x'}^{x''} L_t^x dg_j(x) \to \int_{x'}^{x''}
% L_t^x dg(x),\ \ as\ j\to \infty,
% \end{eqnarray*}
% and
%we get the result.

  %Finally, we deduce an extension of Bouleau-Yor formula by a smoothing procedure. It
%is
%  easy to
%  see we can choose a continuous version of $\int_{-\infty}^{\infty}\nabla ^-f(x)d_x L_t^x$
%  as $-f(X_t)+f(X_0)+\int_0^t \nabla
%^-f(X_s)dX_s$ is a continuous process and
%$$P(\int_{-\infty}^{\infty}\nabla ^-f(x)d_x L_t^x=-f(X_t)+f(X_0)+\int_0^t \nabla
%^-f(X_s)dX_s)=1$$
% for any $t$.
%With the definition of the integral of local time and the
%convergence results of the integrals,
Applying the standard smoothing
procedure on $f(x)$, we can get  $f_n(x)$ which is defined in the same way as
 $g_j(x)$ in (\ref{levynew62}). And by It$\hat {\rm o}$'s formula (c.f. \cite{protter}), we have
 \begin{eqnarray}\label{levynew65}
f_n(X_t)=f_n(X_0)+\int_0^t \nabla
f_n(X_s)dX_s+ A_t^n,\ \ 0\le t<\infty.
\end{eqnarray}
where 
\begin{eqnarray}\label{2ajune2009}
A_t^n={1\over 2}\int_0^t f_n''(X_{s-})d[X,X]_s^c
+\sum_{0\leq s\leq t}[f_n(X_s)-f_n(X_{s-})-\nabla
f_n(X_{s-})\Delta X_s].
\end{eqnarray}
%Let $n\to \infty$ in (\ref{levynew65}), we obtain,
 %\begin{eqnarray}\label{levynew67}
%f(X_t)=f(X_0)+\int_0^t \nabla
%^-f(X_s)dX_s+ A_t,\ \ 0\le t<\infty,
%\end{eqnarray}
%where $\lim\limits_{n\to\infty}A_t^n=A_t$ in $L^2$, because other terms converge in $L^2$. 
From the 
occupation times formula, the definition of the integral of local time and the
convergence results of the integrals, we have 
\begin{eqnarray*}
\lim_{n\to \infty}{1\over 2}\int_0^t f_n''(X_{s-})d[X,X]_s^c
=\lim_{n\to \infty}{1\over 2}\int_{-\infty}^{\infty}L_t^x d\nabla ^-f_n(x)
%=-\lim_{n\to \infty}{1\over 2}\int_{-\infty}^{\infty}\nabla ^-f_n(x)d_x
%L_t^x
=-{1\over 2}\int_{-\infty}^{\infty}\nabla ^-f(x)d_x
L_t^x,
\end{eqnarray*}
and the rough path integral $\int_{-\infty}^{\infty}\nabla ^-f(x)d_x
L_t^x$ is continuous in $t$  from Proposition \ref{may5}. 
For the convergence of jump part in (\ref{2ajune2009}), we can conclude from the proof of
Theorem 3 in Eisenbaum and Kyrianou \cite{eisenbaum12}, if the following assumption:
\medskip

{\it Condition (A):  $\int _{\{|y|<1\}}|f(x+y)-f(x)-\nabla ^-f(x)y|n(dy)$ is well defined and locally bounded in $x$},
\medskip
holds, then 
\begin{eqnarray*}
\lim _{n\to \infty}\sum_{0\leq s\leq t}[f_n(X_s)-f_n(X_{s-})-\nabla
f_n(X_{s-})\Delta X_s]= \sum_{0\leq s\leq t}[f(X_s)-f(X_{s-})-\nabla^-
f(X_{s-})\Delta X_s],
\end{eqnarray*}
in $L^2(dP)$. 
%So we can get the
%\begin{eqnarray}
%A_t=-{1\over 2}\int_{-\infty}^{\infty}\nabla ^-f(x)d_x
%L_t^x
%+\sum_{0\leq s\leq t}[f(X_s)-f(X_{s-})-\nabla
%^-f(X_{s-})\Delta X_s].
%\end{eqnarray}
Therefore we have:
 \begin{thm}\label{theo1}
 Let $f:R\to
 R$ be an absolutely continuous function and have left
 derivative $\nabla^-f(x)$ being left continuous and locally
 bounded, $\nabla^-f(x)$ be of bounded $q$-variation, where $1\leq
 q<3$. Then for $X=(X_t)_{t\geq 0}$, a time homogeneous L\'evy process with $\sigma\neq 0$ and L\'evy measure $n(dy)$
 satisfying Condition (A), and (\ref{levy110}) when $1\leq q <2$, (\ref{levynew36}) when $2\leq
 q<3$, we have P-a.s.
 \begin{eqnarray}\label{may}
f(X_t)&=&f(X_0)+\int_0^t \nabla
^-f(X_s)dX_s-{1\over 2}\int_{-\infty}^{\infty}\nabla ^-f(x)d_x
L_t^x\nonumber\\
&&+\sum_{0\leq s\leq t}[f(X_s)-f(X_{s-})-\nabla
^-f(X_{s-})\Delta X_s], \ \ 0\le t<\infty.
\end{eqnarray}
Here the integral $\int_{-\infty}^{\infty}\nabla ^-f(x)d_x L_t^x$ is
a Lebesgue-Stieltjes integral when $q=1$, a Young integral when
$1<q<2$ and a Lyons' rough path integral when $2\leq q<3$
respectively. In particular, under the condition (\ref{levynew38}),
(\ref{may}) holds for any $2\leq q<3$.
\end{thm}

{\bf Acknowledgements}

It is our great pleasure to thank K.D. Elworthy, W.V. Li, T. Lyons, Z.M. Ma,
S.G. Peng, Z.M. Qian for stimulating conversations. CF would like to
acknowledge the support of National Basic Research Program of China
(973 Program No. 2007CB814903), National Natural Science
Foundation of China (No. 70671069) and the Mathematical Tianyuan
Foundation of China (No. 10826090), Loughborough University and the
London Mathematical Society that enabled her to visit Loughborough
University.

\footnotesize
\begin{thebibliography} {[99]}

%\bibitem{yor1} J. Az\'ema, T. Jeulin, F. Knight and M. Yor, {\em Quelques calculs de compensateurs impliquant l'injectivit\'e de certauns processus croissants}, S\'eminaire de Probabilit\'es XXXII (1998), LNM1686, 316-327.

\bibitem{apple}  D. Applebaum, {\em L\'evy processes and stochastic calculus},
Cambridge Studies in Advanced Mathematics No.93, Cambridge
University Press (2004).


\bibitem{barlow} M. T. Barlow, {\em Necessary and sufficient conditions
for the continuity of local time of L\'evy processes}, Ann. Prob. 16
(1988), 1389-1427.

\bibitem{bou} N. Bouleau and M. Yor, {\em Sur la variation quadratique des temps locaux de certaines semimartingales}, C.R.Acad, Sci. Paris, Ser.I Math 292 (1981), 491-494.

\bibitem{bo} E. S. Boylan, {\em Local times for a class of
Markov processes}, Illinois J. Math., 8: 19-39, 1694.

\bibitem{eisenbaum1} N. Eisenbaum, {\em Integration with respect to local time}, Potential analysis 13
(2000), 303-328.

%\bibitem{eisenbaum2} N. Eisenbaum, {\em Local time-space calculus for revisible semimartingales}, S\'eminaire de Probabilit\'es vol 40, Lectures Notes in Mathematics, Springer-Verlag (to appear).
%

\bibitem{eisenbaum11} N. Eisenbaum, {\em Local time-space stochastic
calculus for L\'evy processes}, Stochastic processes and their
applications 116 (2006), 757-778.

\bibitem{eisenbaum12} N. Eisenbaum and A. Kyprianou, {\em On the parabolic generator of a general one-dimensional L\'evy
process}, Elect. Comm. in Probab. 13 (2008), 198-209.


\bibitem{Zhao1} K. D. Elworthy, A. Truman and H. Z. Zhao, {\em Generalized It$\hat {\rm o}$ Formulae
and space-time Lebesgue-Stieltjes integrals of local times},
S\'eminaire de Probabilit\'es, Vol XL, Lecture Notes in Mathematics
1899, Springer-Verlag, (2007), 117-136.


\bibitem{Zhao33} C. R. Feng and H. Z. Zhao, {\em Two-parameter $p,q$-variation Path and Integration of Local Times},  Potential Analysis,
Vol 25 (2006), 165-204.

\bibitem{Feng3} C.R. Feng and H.Z. Zhao, {\em Rough Path Integral of Local Time},
 C.R. Acad. Sci. Paris, Ser. I 346 (2008), 431-434.

 \bibitem{frw} F. Flandoli, F. Russo and J. Wolf, {\em Some SDEs with distributional drift.
Part II: Lyons-Zheng structure, Ito's formula and semimartingale
characterization}, Random Oper. Stochastic Equations, Vol. 12, No.
2, (2004), 145-184.

\bibitem{Protter} H. F\"ollmer, P. Protter and A. N. Shiryayev, {\em Quadratic covariation and an extension of It$\hat {\rm o}$'s Formula}, Bernoulli 1 (1995), 149-169.

\bibitem{Protter2} H. F\"ollmer and P. Protter, {\em On It$\hat {\rm o}$'s Formula for multidimensional Brownian motion}, Probability Theory and Related Fields 116 (2000), 1-20.

\bibitem{geke} R. K. Getoor and H. Kesten, {\em Continuous of
local times for Markov processes.} Compositio Math., 24: 277-303,
1972.




%\bibitem{terry} T. Lyons and Z. Qian, {\em System Control and Rough Paths}, Clarendon Press Oxford, 2002.
%
%\bibitem{meyer} P. A. Meyer, {\em Un cours sur les int\'egrales stochastiques}, S\'em. Probab 10, Lecture Notes in Math, No. 511, Springer-velay (1976), 245-400.
%
%\bibitem{ks} I. Karatzas and S. E. Shreve, {\em Brownian Motion and
%Stochastic Calculus, Second Edition}, Springer-Verlag: New
%York, 1998.
%
%\bibitem{yor} D. Revuz and M. Yor, {\em Continuous Martingales and
%Brownian Motion, Second Edition}, Springer-Verlag: Berlin,
%Heidelberg, 1994.
%
%\bibitem{williams} David R.E. Williams, {\em Path-wise solutions of stochastic differential equations driven by Levy processes}, Rev. Mat. Iberoam, 17 (2001), 295-329.


%\bibitem{bass} R. F. Bass,  B. M. Hambly and T. J. Lyons (1998),
%{\em Extending the Wong-Zakai theorem to reversible Markov
%processes,} J. Euro. Math. Soc., 4(2002), 237-269.
%
%
%\bibitem{brosamler} G. A. Brosamler, {\em Quadratic variation of
%potentials and harmonic functions,} Transactions of the American
%Mathematical Society 149, 243-257, 1970.
%
%\bibitem{Chung} K. L. Chung, R. J. Williams, {\em Introduction to Stochastic Integration}, Birkhauser
%1990.



%\bibitem{Zhao2} K. D. Elworthy, A. Truman and H. Z. Zhao, {\em
%Asymptotics of Heat Equations with Caustics in One-Dimension},
%Preprint (2006).





%\bibitem{peskir1} R. Ghomrasni and G. Peskir, {\em
%Local time-space calculus and extensions of It$\hat o$'s formula},
%Progr. Probab. Vol. 55, Birkhauser Basel, (2003), 177-192.
%
%\bibitem{hambly} B. M. Hambly any T. L. Lyons, {\em
%Stochastic area for Brownian motion the Sierpinski gasket,} Ann.
%Prob., 26 (1998), 132-48.

\bibitem{ikeda} N. Ikeda and S. Watanabe, {\em Stochastic
Differential Equations and Diffusion Processes}, 2nd Edition,
North-Holland Publ. Co., Amsterdam Oxford New York; Kodansha Ltd.,
Tokyo, 1981.

\bibitem{ito1} K. It$\hat {\rm o}$  {\em Stochastic Integral}, Proc. Imperial Acad. Tokyo 20, 519-524, 1944.

%\bibitem{ito2} K. It$\hat {\rm o}$  {\em On a Stochastic Integral Equation}, Proc. Imperial Acad. Tokyo 22, 32-35, 1946.
%
%\bibitem{ito3} K. It$\hat {\rm o}$  {\em On Stochastic Differential Equations}, Mem. Amer. Mat. Soc. 4, 1-51.
%
%\bibitem{ito4} K. It$\hat {\rm o}$  {\em Lectures on Stochastic Processes}, Tata Institute of Fundamental Research, Vol 24, Bombay, 1961.
%
%\bibitem{ito5} K. It$\hat {\rm o}$ and H. P. Mckean {\em Diffusion Processes and Their Sample Path}, Springer-Verlag, Berlin, 1974.

\bibitem{ks}I. Karatzas and S. E. Shreve, {\em Brownian Motion and
Stochastic Calculus, Second Edition}, Springer-Verlag: New York,
1998.

%\bibitem{kendall} W. Kendall, {\em The radial part of a
%$\Gamma$-martingale and non-implosion theorem}, The Annals of
%Probability, 23, 479-500, 1995.

%\bibitem{kun}H. Kunita, {\em Stochastic Flows and Stochastic Differential Equations}, Cambridge University Press, Cambridge, 1990.

\bibitem{kunwat} H. Kunita and S. Watanabe, {\em On square-integrable
martingales}, Nagoya Math. J., Vol 30 (1967), 209-245.


%\bibitem{ledoux} M. Ledoux, T. J. Lyons and Z. Qian, {\em
%L\'evy area of Wiener processes in Banach Spaces}, Ann. Prob., Vol.
%30, 546-578, 2002.

\bibitem{lejay} A. Lejay, {\em An Introduction to Rough Paths},
S\`eminaire de probabilit\`es XXXVII , Lecture Notes in Mathematics,
Vol. 1832,  1-59, Springer-Verlag, 2003.


\bibitem{levy} P. L\'evy, {\em Processus Stochastiques et Mouvement Brownien}, Gauthier-Villars, Paris, (1948).

\bibitem{terry2} T. Lyons, {\em Differential equations driven by rough signals (I): An extension of an inequality of L. C. Young}, Math. Res. Lett., 1, 451-64, 1994.

\bibitem{terry1} T. Lyons, {\em Differential equations driven by rough signals}, Rev. Mat. Iberoamer., 14, 215-310, 1998.

\bibitem{terry} T. Lyons and Z. Qian, {\em System Control and Rough Paths}, Clarendon Press Oxford, 2002.

\bibitem{rosen1} M. B. Marcus and J. Rosen, {\em Sample path properties of
the local times of strongly symmetric Markov processes via Gaussian
processes}, Ann. Prob. Vol. 20 (1992), 1603-1684.

\bibitem{rosen} M. B. Marcus and J. Rosen, {\em $p$-variation of
the local times of symmetric stable processes and of Gaussian
processes with stationary increments}, Ann. Prob. Vol.20 (1992),
1685-1713.

%\bibitem{mc} E. McShane, {\em Integration}, Princeton University
%Press, Princeton, 1944.

\bibitem{meyer} P. A. Meyer, {\em Un cours sur les int\'egrales stochastiques}, S\'em. Probab 10, Lecture Notes in Math, No. 511, Springer-velay (1976), 245-400.

 \bibitem{nualart} S. Moret and
D. Nualart, {\em Generalization of It$\hat {\rm o}$'s formula for
smooth nondegenerate  martingales}, Stochastic Process. Appl. 91,
115-149, 2001.

\bibitem{perkins} E. Perkins, {\em Local time is a semimartingale}, Z.
Wahrsch. Verw. Gebiete, Vol. 60 (1982), 79-117.

\bibitem{protter} P. E. Protter, {\em Stochastic Integration and Differential Equations}, 2nd edition, Springer, 2005. 

%\bibitem{peskir2} G. Peskir, {\em A change-of-variable formula with local time on curves}, J.Theoret
%Probab. (to appear).

%\bibitem{peskir3} G. Peskir, {\em A change-of-variable formula with local time on surfaces},
%S\'eminaire de Probabilit\'es vol 40, Lectures Notes in Mathematics,
%Springer-Verlag (to appear).

%\bibitem{peskir4} G. Peskir, {\em On the American option problem}, Math.Finance, vol 15 (2005), 169-181.

\bibitem{yor} D. Revuz and M. Yor, {\em Continuous Martingales and
Brownian Motion, Second Edition}, (Springer-Verlag: Berlin,
Heidelberg, 1994).

\bibitem{rog} L. C. G. Rogers and J. B. Walsh, {\em Local time and stochastic area integrals}, Annals of Probas. 19(2) (1991), 457-482.

%\bibitem{rw} L. C. G. Rogers and D. Williams, {\em Diffusions, Markov Processes and Martingales}, Vol. 2,
%It$\hat {\rm o}$ Calculus, Cambridge University Press, 2nd edition,
%2000.
%
%\bibitem{rv} F. Russo and P. Vallois, {\em It$\hat {\rm o}$'s Formula
%for $C^1$-functions of semimartingales}, Probability Theory and
%Related Fields, 104 (1996), 27-41.

%\bibitem{ste} E. M. Stein, {\em Singular Integrals and Differentiability Properties of Functions}, Princeton
%University Press, 1970.

\bibitem{tan} H. Tanaka, {\em Note on continuous additive functionals of the 1-dimensional Brownian path}, Z.Wahrscheinlichkeitstheorie and Verw Gebiete 1 (1963), 251-257.

%\bibitem{walsh} J. B. Walsh, {\em An Introduction to Stochastic
%Partial Diffenrential Equations}, In \'Ecole d\'et\'e de
%Probabilit\'e de Saint Flour, XIV, ed. D.L. Hennequin, Lecture Notes
%in Mathematics No. 1180, (1986), 265-439.

\bibitem{wang} A. T. Wang, {\em Generalized It$\hat {\rm o}$'s formula and additive functionals of Brownian motion}, Z.Wahrscheinlichkeitstheorie and Verw Gebiete, 41(1977), 153-159.

\bibitem{williams} D. R. E. Williams, {\em Path-wise solutions of stochastic defferential equations driven by Levy processes}, Rev. Mat. Iberoam, 17 (2001),
295-329.

\bibitem{young1} L. C. Young, {\em An inequality of H${\rm \ddot{o}}$lder type, connected with Stieltjes integration}, Acta Math., 67 (1936), 251-282.



\end{thebibliography}
\end{document}